\documentclass[aop,MSNbibl,seceqn,citesort,dvips]{arximspdf}

\doi{10.1214/11-AOP655}
\volume{40}
\issue{4}
\pubyear{2012}
\firstpage{1436}
\lastpage{1482}

\makeatletter
\let\epsilon\varepsilon

\newtheorem{theorem}{Theorem}[section]
\newtheorem{corollary}[theorem]{Corollary}
\newtheorem{lemma}[theorem]{Lemma}

\def\R{\mathbb{R}}
\def\Z{\mathbb{Z}}
\def\E{\mathbb{E}}
\def\ol{\overline}
\def\P{\mathbb{P}}

\makeatother

\begin{document}
\begin{frontmatter}

\title{Quenched asymptotics for Brownian motion of renormalized
Poisson potential and for the related parabolic Anderson models}
\runtitle{Brownian motion of renormalized Poisson potential}

\begin{aug}
\author[A]{\fnms{Xia} \snm{Chen}\corref{}\thanksref{aut1}\ead[label=e1]{xchen@math.utk.edu}}
\thankstext{aut1}{Supported in part by NSF Grant DMS-07-04024.}
\runauthor{X. Chen}
\affiliation{University of Tennessee}
\address[A]{Department of Mathematics\\
University of Tennessee\\
Knoxville, Tennessee 37996\\
USA\\
\printead{e1}} 
\end{aug}

\received{\smonth{9} \syear{2010}}
\revised{\smonth{1} \syear{2011}}

%
\begin{abstract}
Let $B_s$ be a $d$-dimensional Brownian motion and $\omega(dx)$
be an independent Poisson field on $\mathbb{R}^d$.
The almost sure asymptotics for the logarithmic moment
generating function
\[
\log\mathbb{E}_0\exp\biggl\{\pm\theta\int_0^t\overline{V}(B_s)\,ds\biggr
\}\qquad (t\to\infty)
\]
are investigated in connection with
the renormalized Poisson potential of the form
\[
\overline{V}(x)=\int_{\mathbb{R}^d}{1\over\vert
y-x\vert^p}[\omega(dy)-dy], \qquad x\in\mathbb{R}^d.
\]
The investigation is motivated by some practical
problems arising from the models of Brownian motion in
random media and from the parabolic Anderson models.
\end{abstract}

%

\end{frontmatter}

\section{Introduction} \label{intro}

Consider a particle doing a random movement in the space~$\R^d$. The trajectory
of the particle is described by a $d$-dimensional Brownian motion~$B_s$.
Independently, there is a family of the obstacles
randomly located in the space~$\R^d$. Assume that each obstacle has
mass 1 and that the obstacles are distributed in $\R^d$ according to
a Poisson field $\omega(dx)$ with the Lebesgue
measure $dx$ as its intensity measure. Throughout,
the notation ``$\P$'' and ``$\E$'' are used
for the
probability law and the expectation, respectively, generated by the Poisson
field $\omega(dx)$, while the notation ``$\P_x$'' and ``$\E_x$'' are
for the probability law and the expectation, respectively, of the Brownian
motion $B_s$ with $B_0=x$.

The model of Brownian motion in Poisson potential has been
introduced to describe the trajectory of a Brownian particle that survived
being trapped by the obstacles.\vadjust{\goodbreak} We refer the reader to
the book by Sznitman \cite{Sznitman} and the survey \cite{Komorowski}
made by Komorowski for a systematic account of this model and
the monograph by Harvlin and Ben Avraham \cite{HA} for physicists' views
on the trapping kinetics.
In the usual set-up, the random field (known as potential function)
%
\begin{equation}\label{intro-1}
V(x)=\int_{\R^d}K(y-x)\omega(dy)
\end{equation}
represents the total trapping energy
at $x\in\R^d$ generated by the Poisson obstacles, where
$K(x)\ge0$ is a deterministic function on $\R^d$ known as
the shape function.
In the quenched setting, where the observation of the system is
conditioned on the environment generated by the Poisson obstacles,
the model of Brownian motion in Poisson potential is often
introduced as the random Gibbs measure $\mu_{t,\omega}$ on
$C \{[0, t];\R^d \}$ defined as
%
\begin{equation}\label{intro-2}
{d\mu_{t,\omega}\over d\P_0}={1\over Z_{t,\omega}}\exp\biggl\{-\theta
\int_0^tV(B_s)\,ds \biggr\}.
\end{equation}

The integral
\[
\int_0^tV(B_s)\,ds
\]
measures the total trapping energy received by the Brownian particle
up to the time $t$. Under the law $\mu_{t,\omega}$,
therefore, the Brownian paths heavily impacted by the Poisson obstacles
are penalized and become less likely.

Sznitman \cite{Sznitman}
considers two kinds of shape functions. In one case
$K(x)=\infty1_C$ for a nonpolar set $C\subset\R^d$, while in another
case, the shape function~$K(x)$ is assumed to be bounded and compactly supported.
The correspondent
potential functions are called hard and soft obstacles, respectively.
In the case of hard obstacles, the Brownian particle is completely free
from the influence of the obstacles until hitting the $C$-neighborhood
of the
Poisson cloud which serves as the death trap. In the setting of the soft
obstacles, only the obstacles in a local neighborhood of the Brownian
particle act on the particle, and the collision does not create extreme
impact.

According to
Newton's law of universal attraction, for example,
the integrals
\[
\int_{\R^d}{1\over\vert y-x\vert^{d-1}}\omega(dy)\quad\mbox{and}\quad
\int_{\R^d}{1\over\vert y-x\vert^{d-2}}\omega(dy),\qquad
x\in\R^d,
\]
represent (up to constant multiples), respectively,
the total gravitational force and the total gravitational
potential at the location $x$ in the gravitational field generated
by the Poisson obstacles in the case when $d\ge3$.
Therefore, it makes
sense in physics to consider the shape function of the
form\looseness=-1
%
\begin{equation}\label{intro-4}
K(x)=\vert x\vert^{-p},\qquad  x\in\R^d.
\end{equation}\looseness=0

A serious problem is that under choice (\ref{intro-4}),
$V(x)$ blows up at every $x\in\R^d$ when $p\le d$. In a recent paper
\cite{C-K}, a renormalized model has\vadjust{\goodbreak} been proposed as follows:
First, it is shown (\cite{C-K}, Corollary 1.3) that under the assumption
$d/2<p<d$
the renormalized potential
%
\begin{equation}\label{intro-6}
\ol{V}(x)=\int_{\R^d}{1\over\vert y-x\vert^p} [\omega(dy)-dy ],\qquad
x\in\R^d,
\end{equation}
can be properly defined and that for any $\theta>0$ and $t>0$,
\[
\E\otimes\E_0\exp\biggl\{-\theta\int_0^t\ol{V}(B_s)\,dx \biggr\}<\infty.
\]
Consequently,
%
\begin{equation}\label{intro-7}
\ol{Z}_{t,\omega}\equiv
\E_0\exp\biggl\{-\theta\int_0^t\ol{V}(B_s)\,dx \biggr\}<\infty\qquad\mbox{a.s.}
\end{equation}
Thus, the Gibbs measure $\ol{\mu}_{t,\omega}$ given as
%
\begin{equation}\label{intro-7'}
{d\ol{\mu}_{t,\omega}\over d\P_0}={1\over\ol{Z}_{t,\omega}}
\exp\biggl\{-\theta\int_0^t\ol{V}(B_s)\,ds \biggr\}
\end{equation}
is well defined and appears to be a natural extension of $\mu_{t,\omega}$
[given in (\ref{intro-2})] in the following sense: When $K(x)$ is compactly
supported and bounded, by translation invariance of Lebesgue measure,
\begin{eqnarray*}
\ol{V}(x)&=&\int_{\R^d}K(y-x) [\omega(dy)- dy ]
=\int_{\R^d}K(y-x)\omega(dy)-\int_{\R^d}K(y)\,dy\\
&=&V(x)-\mbox{constant}.
\end{eqnarray*}
So the Gibbs measures generated by $V(x)$ and by $\ol{V}(x)$ are equal.
We call the random path under the law $\mu_{t,\omega}$ the Brownian motion
of the renormalized Poisson potential $\ol{V}(x)$.
In the case when $K(x)$ is given
in (\ref{intro-4}), the renormalized Poisson potential
$\ol{V}(x)$ in (\ref{intro-6}) appears as
the constant multiple of the Riesz potential of the compensated
Poisson field $\omega(dy)- dy$.

One of major objectives of this paper is to investigate the large-$t$
asymptotics for partition function $\ol{Z}_{t,\omega}$ given in (\ref{intro-7})
with the potential function~$\ol{V}(x)$ be defined in (\ref{intro-6}).

This problem is also motivated by the parabolic Anderson formulated in
the form
of the Cauchy problem
%
\begin{equation} \label{intro-8}
\cases{
\partial_tu(t,x) =\kappa\Delta u(t,x)+\xi(x)u(t,x),\cr
u(0, x)=1,
}
\end{equation}
where $\kappa>0$ is a constant called diffusion
coefficient, and $\xi(x)$ is a properly chosen random field
called potential.

Among other things,
the parabolic Anderson models are used to describe evolution of the mass
density $u(t,x)$ distributed in $\R^d$ (see, e.g., \cite{C-K} for the
discussion on this link). The mathematical relevance of
the parabolic Anderson\vadjust{\goodbreak} models to our topic is based on two facts: First,
by the space homogeneity of the Poisson field,
%
\begin{equation}\label{intro-15}
\{\xi(t, x); t\ge0 \}\stackrel{d}{=} \{\xi(t, 0); t\ge0 \},\qquad
x\in\R^d.
\end{equation}
Consequently, the focus of the investigation is often on $u(t,0)$. Second,
by the
Feynman--Kac representation,
%
\begin{equation}\label{intro-10}\qquad
u(t,0)=\E_0\exp\biggl\{\int_0^{t}\xi(B_{2\kappa s})\,ds \biggr\}
=\E_0\exp\biggl\{(2\kappa)^{-1}\int_0^{2\kappa t}\xi(B_s)\,ds \biggr\}
\end{equation}
for sufficiently nice $\xi(x)$.

There are long lists of publications 
on this model among which we refer the reader to the monograph \cite{CM}
by Carmona and Molchanov for the overview and background of this subject.
In the usual set-up, $\xi(x)=\pm V(x)$ with~$V(x)$
being given in (\ref{intro-1}). In the existing literature,
the shape function~$K(x)$ is usually assumed to be bounded and compactly
supported so that the potential function $V(x)$ can be
defined. A localized shape is
analogous to the usual set-up in the discrete parabolic Anderson model, where
the potential $\{V(x); x\in\Z^d\}$ is an i.i.d. sequence.
On the other hand, there are practical
needs for considering the cases, such that when $K(x)=\vert x\vert^{-p}$,
where the environment has a long-range dependency and the extreme
force surges at the locations of the Poisson obstacles.

In this paper, we consider the case when $\xi(x)=\pm\theta\ol{V}(x)$
where $\ol{V}(x)$ is defined in (\ref{intro-6}).
Given the fact (Proposition 2.7 in \cite{C-K})
that $\ol{V}(x)$ is
unbounded in any neighborhood with positive probability, it is unlikely
that equation~(\ref{intro-8}) is solvable in the path-wise sense.
On the other hand, it has been proved in~\cite{C-K} that $u(t,x)$ represented
by the Feynman--Kac formula is a~mild solution to~(\ref{intro-8}) [with
$\xi(x)=\pm\theta\ol{V}(x)$]
whenever the quenched moment in~(\ref{intro-10}) is finite.

The objects of our investigation are the quenched
exponential moments
%
\begin{equation}\label{intro-21}
\E_0\exp\biggl\{-\theta\int_0^t\ol{V}(B_s)\,ds \biggr\} \quad\mbox{and}\quad
\E_0\exp\biggl\{\theta\int_0^t\ol{V}(B_s)\, ds \biggr\}.
\end{equation}
According to (\ref{intro-7}), the first exponential moment in (\ref{intro-21})
is almost
surely defined. As for the second exponential moment,
it has been proved in recent work
\cite{C-K} that
the correspondent annealed exponential moment blows up, and that,
for any $\theta>0$ and $t>0$,
%
\begin{equation} \label{intro-22}
\E_0\exp\biggl\{\theta\int_0^t\ol{V}(B_s)\,ds \biggr\}
\cases{
<\infty, &\quad if $p<2$,\cr
=\infty,&\quad if $p>2$.
}
\end{equation}
The critical case $p=2$, in which $d=3$ by the constraint $d/2<p<d$,
has been investigated in a more recent paper \cite{CR} where it is shown
that for\break any $t>0$
%
\begin{equation} \label{intro-23}
\E_0\exp\biggl\{\theta\int_0^t\ol{V}(B_s)\,ds \biggr\}
\cases{
<\infty,&\quad  a.s.
 when $ \theta<{1\over16}$,\vspace*{2pt}\cr
=\infty,&\quad a.s.
when $  \theta> {1\over16}$.
}\vadjust{\goodbreak}
\end{equation}

The main objective of this paper is to investigate the quenched
large-$t$ asymptotics for the exponential moments given in
(\ref{intro-21}) whenever these moments are finite,
except the critical case described in (\ref{intro-23}) (which is
studied in
\cite{CR}). We point out the references
\cite{G-M,DM,CVM,CTV,CM,CV,CMS,Englander,FV,G-K-M,HKS,Sznitman-1,Sznitman-2,Sznitman,VZ}
as an incomplete list related to this topic.

For later comparison, we mention some existing results which are
narrowly relevant to the topic of this paper. Let
the potential function $V(x)$ be given in (\ref{intro-1}).
Sznitman (\cite{Sznitman}, Theorem 5.3,
page 196) shows that for the bounded and compactly supported
shape $K(\cdot)$ and $\theta>0$,
%
\begin{equation}\label{intro-18}
\qquad \lim_{t\to\infty}{(\log t)^{2/d}\over t}\log\E_0\exp\biggl\{-\theta\int_0^t
V(B_s)\,ds \biggr\}
=-\lambda_d \biggl({\omega_d\over d} \biggr)^{2/d}\quad  \mbox{a.s.-}\P,\hspace*{-20pt}
\end{equation}
where $\lambda_d>0$ is
the principal eigenvalue of the Laplacian operator
$(1/2)\Delta$ on the $d$-dimensional unit ball
with zero boundary values, and $\omega_d$ is the volume of
the $d$-dimensional unit ball. With a slightly different formulation
\cite{Sznitman-1},
the model of hard obstacles yields the same pattern of asymptotics.

Under some continuity, boundedness assumptions on $K(x)$
and under some restriction on the tail of $K(x)$,
Carmona and Molchanov (\cite{CM-1}, Theorem~5.1) prove that
%
\begin{equation}\label{intro-19}
\qquad \lim_{t\to\infty}{\log\log t\over t\log t}
\log\E_0\exp\biggl\{\theta\int_0^tV(B_s)\,ds \biggr\}
=d\theta\sup_{x\in\R^d}K(x)\quad  \mbox{a.s.-}\P.
\end{equation}
The interested reader is also referred to
\cite{G-M-1} and \cite{G-K-M} for the correspondent
asymtotics of the second order.

After the first draft of this paper was completed,
the author learned the recent investigation
by Fukushima \cite{Fu} in the case
when \mbox{$K(x)=\vert x\vert^{-p}\wedge1$}\break with $d<p<d+2$, the setting where no
renormalization is necessary. Fukushi\-ma~\cite{Fu} shows that
\begin{eqnarray} \label{intro-20}
&&\lim_{t\to\infty}t^{-1}(\log t)^{-{(p-d)/d}}\log\E_0\exp\biggl\{-\int_0^t
V(B_s)\,ds \biggr\}\nonumber\\[-8pt]\\[-8pt]
&&\qquad=-{d\over p} \biggl({p-d\over pd} \biggr)^{(p-d)/d}
\biggl(\omega_d\Gamma\biggl({p-d\over d} \biggr) \biggr)^{p/d}\qquad  \mbox{a.s.-}\P.\nonumber
\end{eqnarray}
It should be mentioned that Fukushima also obtained the second
asymptotic term in his setting.

\section{Main theorems and strategies}\label{theorem}

Throughout this paper, let $\omega_d$ be the volume of the $d$-dimensional
unit ball. Let $W^{1,2}(\R^d)$ denote the Sobolev space
given as
\[
W^{1,2}(\R^d)= \{f\in\mathcal{L}^2(\R^d); \nabla f\in
\mathcal{L}^2(\R^d) \}.
\]
By (\ref{theorem-6}) below, when $d/2<p<\min\{d, 2\}$
there is a constant $C>0$ such that
\[
\int_{\R^d}{f^2(x)\over\vert x\vert^p}\,dx\le C\|f\|_2^{2-p}\|\nabla
f\|_2^p,
\qquad f\in W^{1,2}(\R^d).
\]
Let $\sigma(d,p)>0$ be the best constant in above inequality.

The main theorems are stated as follows.

\begin{theorem}\label{th-2} Under $d/2<p<d$,
\begin{eqnarray} \label{theorem-10}
&&\lim_{t\to\infty}t^{-1}(\log t)^{-{(d-p)/ d}}\log\E_0\exp\biggl\{-\theta
\int_0^t
\ol{V}(B_s)\,ds \biggr\}\nonumber\\[-8pt]\\[-8pt]
&&\qquad={\theta d^2\over d-p}
\biggl({\omega_d\over d}\Gamma\biggl({2p-d\over p} \biggr) \biggr)^{p/d}\qquad
\mbox{a.s.-}\P\nonumber
\end{eqnarray}
for every $\theta>0$.
\end{theorem}

\begin{theorem}\label{th-3} Under $d/2<p<\min\{2, d\}$,
\begin{eqnarray} \label{theorem-11}
\qquad &&\lim_{t\to\infty}{1\over t} \biggl({\log\log t\over\log t} \biggr)^{2/(2-p)}
\log\E_0\exp\biggl\{\theta\int_0^t\ol{V}(B_s)\,ds
\biggr\}\nonumber\\[-8pt]\\[-8pt]
\qquad &&\qquad={1\over2}p^{p/(2-p)}
(2-p)^{(4-p)/(2-p)} \biggl({d\theta\,\sigma(d,p)\over2+d-p} \biggr
)^{2/(2-p)}\qquad
\mbox{a.s.-}\P\nonumber
\end{eqnarray}
for every $\theta>0$.
\end{theorem}

We now make a comparison of ``(\ref{intro-18}) versus (\ref
{theorem-10})'' and
``(\ref{intro-19}) versus (\ref{theorem-11}).''
First, the quenched exponential
moments in
our models generate significantly larger quantities. Second, a heavy
shape dependence (or $p$-dependence) presented in
our theorems sharply contrasts the shape
insensitivity appearing in (\ref{intro-18}) and (\ref{intro-19}).
In Theorem \ref{th-2},
it is the nonlocality of the shape function that plays a major role, while
the high peaks of $\ol{V}(x)$ correspond to small values of
the quenched exponential moment.
On the other hand, the asymptotics in Theorem \ref{th-3}
is shaped by the singularity of $K(x)=\vert x\vert^{-p}$ at $x=0$.
In addition, there seems
to be a degree of resemblance between (\ref{intro-20}) and (\ref{theorem-10}).
Based on the comment made about roles of nonlocality and singularity,
it may be possible that (\ref{intro-20}) remains true even without removing
the singularity of $K(x)$ at $x=0$. We leave this problem to future study.

Does the Lebesgue measure in renormalization contribute to the
limit laws stated in Theorems \ref{th-2} and  \ref{th-3}?
The answer is ``Yes'' to Theorem \ref{th-2}, for otherwise the
right-hand side
of (\ref{theorem-10}) would be negative. The answer is ``No'' to
Theorem \ref{th-3} as the major impact comes from the Poisson points
in a very small neighborhood
of the site where the Brownian particle is located [see (\ref
{theorem-22}) below
for a more quantified analysis on this point].\vadjust{\goodbreak}

Associated with the spatial Brownian motion in
the classic gravitational field generated by
the Poisson obstacles, the following corollary appears as Theorem~\ref{th-2}
in the special case $d=3$ and $p=2$.\vspace*{-3pt}

\begin{corollary}\label{th-6}
When $d=3$ and $p=2$,
%
\begin{equation} \label{theorem-20}
\qquad \lim_{t\to\infty}t^{-1}(\log t)^{-1/3}\log\E_0\exp\biggl\{-\theta\int_0^t
\ol{V}(B_s)\,ds \biggr\}=3\root3\of{12}\pi\theta
\quad\mbox{a.s.-}\P
\end{equation}
for every $\theta>0$.\vspace*{-3pt}
\end{corollary}

Let $u_0(t,x)$ and $u_1(t,x)$ be the mild solutions to the parabolic Anderson\vspace*{1pt}
problems (\ref{intro-8}) that satisfy
the Feynman--Kac representation (\ref{intro-10})
with $\xi(x)=-\theta\ol{V}(x)$ and $\xi(x)=\theta\ol{V}(x)$, respectively.
By the space homogeneity~({\ref{intro-15}) and by Theorems~\ref{th-2}
and \ref{th-3},
\begin{eqnarray} \label{theorem-13}
&&\lim_{t\to\infty}t^{-1}(\log t)^{-{(d-p)/d}}\log u_0(t,x)\nonumber\\[-9pt]\\[-9pt]
&&\qquad ={\theta d^2\over d-p}
\biggl({\omega_d\over d}\Gamma\biggl({2p-d\over p} \biggr) \biggr)^{p/d}\quad
\mbox{a.s.-}\P,\nonumber\vspace*{-3pt}
\end{eqnarray}
\begin{eqnarray}
\label{theorem-14}
\qquad &&\lim_{t\to\infty}{1\over t} \biggl({\log\log t\over\log t} \biggr
)^{2/(2-p)}
\log u_1(t,x)\nonumber\\[-9pt]\\[-9pt]
\qquad &&\qquad ={1\over2} \biggl({p\over2\kappa} \biggr)^{p/(2-p)}
(2-p)^{(4-p)/(2-p)} \biggl({d\theta\,\sigma(d,p)\over2+d-p} \biggr
)^{2/(2-p)}\quad
\mbox{a.s.-}\P\nonumber
\end{eqnarray}
for every $\theta>0$ and $x\in\R^d$.

An immediate observation is that the diffusion coefficient $\kappa$ does
not appear in (\ref{theorem-13}). The same phenomena have been noticed
by Carmona and Molchanov \cite{CM-1}
in the case when $\xi(x)=\theta V(x)$ for the same $V(x)$ appearing
in (\ref{intro-19}).

In the following we compare the strategies for the laws given
in (\ref{intro-18}), (\ref{intro-19}), (\ref{theorem-10})
and (\ref{theorem-11}). To make the discussion more informative, we
focus on
the lower bounds and try to describe the behavior of the Brownian
particle and the behavior of the Poisson particle in each strategy.
The treatment for~(\ref{intro-18}) and~(\ref{intro-19}) does not have
to be
the same as their original proof. In our discussion, we use the notation
$B(x, R)$ for the $d$-dimensional ball of the center $x$ and radius $R$.

\setcounter{footnote}{1}
The following ingredients on the behavior of the Brownian particle are
common to all strategies:
Up to the time $t$ the Brownian particle
stays in the ball $B(0, R_t)$ (referred as
``macro-ball'') with the radius $R_t$ roughly equal to
$t$.\footnote{The combination of the word ``roughly'' and a big number $t$
means $tL(t)$ with $L(t)$ slow-varying at $\infty$.}\vadjust{\goodbreak} Within a period
$[0,ut]$ (with a very small $u>0$), the Brownian particle moves into
one of the roughly $t^d$ prearranged and evenly located
identical micro-balls
%
\begin{equation}\label{theorem-17'}
D_z\equiv B(z,r_t);\qquad z\in b_t\Z^d\cap B(0, R_t),
\end{equation}
where $r_t\ll b_t$ and $r_tR_t\ll t$. The principle that Brownian
particle chooses~$D_z$ is to
maximize the positive energy
(or to minimize the negative energy) from the Poisson field.

The main difference among different strategies
in the Brownian path is on the radius $r_t$ of
the microbes.
By the relation $r_tR_t\ll t$ and by
a classic small ball
estimate, the cost for the Brownian particle to choose $D_z$ is ($\delta>0$
is a small number here)
%
\begin{eqnarray}\label{theorem-17}
&&\P_0 \{\mbox{The Brownian particle reaches $D_z$ quickly}\nonumber\hspace*{-35pt}\\
&&\hphantom{\P_0\{The Brownian particle reaches}  \mbox{and then stays in $D_z$ up to $t$}\}
\nonumber\hspace*{-35pt}\\[-8pt]\\[-8pt]
&&\quad\!\ge\!{1\over(2\pi)^d}\int_{B(z,\delta r_t)}\!e^{-\vert x\vert^2/(2ut)}
  \P_0\{ B_s\!\in\! B(z\!-\!x, r_t) \mbox{ for $0\!\le\! s\!\le\!(1\!-\!u)t$} \}\,dx\hspace*{-35pt}
  \nonumber\\
&&\quad\!\approx\!\exp\{-o (r_t^{-2}t ) \}
  \P_0 \Bigl\{\sup_{0\le s\le t}\vert B_s\vert\!\le\! r_t \Bigr\}
  \!\approx\!\exp\{-\lambda_d r_t^{-2}t \}.\nonumber\hspace*{-35pt}
\end{eqnarray}
Here we recall that $\lambda_d>0$ is the principle eigenvalue of the Laplacian
operator $(1/2)\Delta$ on the $d$-dimensional unit ball with zero boundary
condition. To make the cost affordable compared with the deviation
scale $t(\log t)^{-2/d}$ in the strategy for (\ref{intro-18}),
for example, the radius
$r_t$ should be at least $r(\log t)^{1/d}$ with the constant $r>0$.
Based on the same principle, the critical
radius of the micro-balls in each strategy are determined as following:
%
\begin{equation}\label{theorem-18}
r_t= \cases{
\displaystyle r (\log t)^{1/d},&\quad  in the strategy for (\ref{intro-18}),\vspace*{2pt}\cr
\displaystyle r\sqrt{\log\log t\over\log t},&\quad
in the strategy for (\ref{intro-19}),\vspace*{2pt}\cr
\displaystyle r (\log t)^{-{(d-p)/(2d)}},&\quad
in the strategy for (\ref{theorem-10}),\vspace*{2pt}\cr
\displaystyle r \biggl({\log\log t\over\log t} \biggr)^{1/(2-p)},&
\quad in the strategy for (\ref{theorem-11}).
}
\end{equation}

We now describe the behavior of the Poisson field in each strategy.
For~(\ref{intro-18}), the high peak of the quenched moment occurs
when $\int_0^tV(B_s)\,ds\approx0$. To make this happen,
one of the $C$-neighborhoods $\widetilde{D}_z\equiv D_z+C$
[$z\in b_t\Z^d\cap B(0, R_t)$] is obstacle-free,
where $C\subset\R^d$ is the compact support of $K(x)$,
and the Brownian particle spends most of its time in that same
micro-ball~$D_z$.
In view of (\ref{theorem-17}), therefore,
\[
\E_0\exp\biggl\{-\theta\int_0^t V(B_s)\,ds \biggr\}
\succeq\exp\{-\lambda_dr^{-2}t(\log t)^{-2/d} \}
\]
on the event $ \{\min_z\omega(\widetilde{D}_z)=0\}$,
where the relation ``$\succeq$'' reads as ``asymptotically greater than
or equivalent to.''

On the other hand,
\begin{eqnarray*}
\P\Bigl\{\min_z\omega(\widetilde{D}_z)=0 \Bigr\}
&\approx&1- \bigl(1-\P\{\omega(\widetilde{D}_0)=0 \} \bigr)^{t^d}\\
&=&1- (1-\exp\{-\omega_dr^d\log t \} )^{t^d}\\
&\approx&1-\exp\{-t^{d-\omega_dr^d} \}.
\end{eqnarray*}
Hence, a standard way of using the Borel--Cantelli lemma
shows that the phase transition
between
\begin{eqnarray}\label{theorem-19}
\P\Bigl\{\min_z\omega(\widetilde{D}_z)=0 \mbox{ eventually} \Bigr\}&=&1
\quad\mbox{and}\nonumber\\[-8pt]\\[-8pt]
\P\Bigl\{\min_z\omega(\widetilde{D}_z)\ge1 \mbox{ eventually} \Bigr\}&=&1\nonumber
\end{eqnarray}
occurs when $r$ satisfies $\omega_dr^d=d$. Consequently, this strategy
leads to the lower
bound requested by (\ref{intro-18}).

In the strategy for (\ref{theorem-10}), only the impact
of the Poisson obstacles within
the distance $a(\log t)^{1/d}$ from the Brownian particle is counted.
To determine constant $a>0$, a crucial problem is whether or not
the high peak can be captured by the ``empty ball'' strategy which means
to make the ball
$B (B_s, a(\log t)^{1/d} )$ [$\approx B (z,a(\log t)^{1/d} )$
as the Brownian particle
stays in $D_z$] free of the Poisson obstacles. Under the ``empty-ball''
strategy,
\begin{eqnarray*}
\int_0^t\ol{V}(B_s)\,ds&\approx&\int_0^tV_1(B_s)\,ds-
t\int_{\{\vert x\vert\le a(\log t)^{1/d}\}}{1\over\vert x\vert^p}\,dx\\
&\approx&-{a^{d-p}\omega_d\over d-p}t(\log t)^{(d-p)/d},
\end{eqnarray*}
where
\[
V_1(x)=\int_{\vert y-x\vert\le a(\log t)^{1/d}\}}
{\omega(dy)\over\vert y-x\vert^p}.
\]

On the other hand, the estimate given in (\ref{theorem-19}) shows that
the largest radius $R$ for one of the balls $B(z, R)$
[$z\in b_t\Z^d\cap B(0, R_t)$] to be obstacle-free is $R=
(\omega_d^{-1}d)^{1/d}(\log t)^{1/d}$. By making $r>0$ sufficiently large
in (\ref{theorem-17}), the best lower bound
that the ``empty-ball'' strategy can offer is
\begin{eqnarray*}
&&\liminf_{t\to\infty}t^{-1}(\log t)^{-{(d-p)/d}}\log
\E_0\exp\biggl\{-\theta\int_0^t
\ol{V}(B_s)\,ds \biggr\}\\
&&\qquad \ge{\theta\over d-p}d^{(d-p)/d}\omega_d^{p/d}\qquad
\mbox{a.s.}
\end{eqnarray*}
under the optimal choice $a=(\omega_d^{-1}d)^{1/d}$.
In comparison with (\ref{theorem-10}), this bound
gives the right rate but not the right constant.

Based on the above analysis, we conclude
that the constant $a>0$ has to be arbitrarily large
and that the ``empty-ball'' strategy is not working
well for (\ref{theorem-10}).

We now come to (\ref{intro-19}). By the continuity assumption on the
shape function and by homogeneity of the Poisson field,
the supremum $ \sup_{x\in\R^d}K(x)$ can be achieved
somewhere, and we may assume that
$K(0)= \sup_{x\in\R^d}K(x)$ in the following discussion.
To support the limit law given in (\ref{intro-19}), the Poisson field executes
a strategy that fills one of the $\delta$-balls $ \{B(z,\delta);
z\in b_t\Z^d\cap B(0, R_t) \}$ with a high density of the Poisson points,
where the constant $\delta>0$ is (arbitrarily) small but fixed.
By
translation invariance and by continuity of $K(x)$, for any
$z\in b_t\Z^d\cap B(0, R_t)$
\begin{eqnarray*}
V(B_s)&=&\int_{\R^d}K(x-B_s)\omega(dx)\\
&=&\int_{\R^d}K \bigl(x-(B_s-z) \bigr)\omega(z+dx)\succeq K(0) \omega
(B(z,\delta) )
\end{eqnarray*}
as $B_s\in D_z$. By (\ref{theorem-17}), therefore,
\begin{eqnarray}\label{theorem-20'}
&&\E_0\exp\biggl\{\theta\int_0^tV(B_s)\,ds
\biggr\}\nonumber\\[-8pt]\\[-8pt]
&&\qquad \succeq\exp\biggl\{\theta K(0) t\max_z\omega(B(z,\delta) )
-\lambda_d r^{-2}{t\log t\over\log\log t} \biggr\}.\nonumber
\end{eqnarray}

On the other hand, by independence
\begin{eqnarray*}
\P\biggl\{\max_z\omega(B(z,\delta) )\ge\sigma{\log t\over\log\log t}
\biggr\}&\approx&1- \biggl(1-\P\biggl\{\omega(B(0,\delta) )
\ge\sigma{\log t\over\log\log t} \biggr\} \biggr)^{t^d}\\
&\approx&1- (1-\exp\{-\sigma\log t \} )^{t^d}\\
&\approx&
1-\exp\{-t^{d-\sigma} \}\qquad \forall\sigma>0.
\end{eqnarray*}
Using the Borel--Cantelli lemma we can prove that
%
\begin{equation}\label{theorem-21}
\lim_{t\to\infty}{\log\log t\over\log t}\max_z\omega(B(z,\delta) )
=d \qquad\mbox{a.s.}
\end{equation}
Since $r>0$ can be arbitrarily large, (\ref{theorem-20'}) and (\ref{theorem-21})
lead to the lower bound requested by (\ref{intro-19}).

The strategy that Poisson field executes in
(\ref{theorem-11})
is to fill one of the balls
\[
B \biggl(z, \delta\biggl({\log\log t\over\log t}
\biggr)^{1/(2-p)} \biggr);\qquad
z\in b_t\Z^d\cap B(0, R_t),
\]
with a high concentration of the Poisson points.
In the following we present a
simple algorithm to illustrate the idea.
Assume that the Brownian particle spends most of its
time in $D_z$ for some $z\in b_t\Z^d\cap B(0, R_t)$.
Given a fixed $a>0$, it is not hard to show that
the impact of
the Poisson points which are $a$-unit away from the Brownian
particle is negligible, and that the ``renormalizer'' does not
make any noticeable contribution to the limit law in (\ref
{theorem-11}). Hence,
%
\begin{eqnarray}\label{theorem-22}
\qquad \ol{V}(B_s)&=&\int_{\R^d}{1\over\vert x-(B_s-z)\vert^p}
[\omega(z+dx)-dx ]\nonumber\\[-2pt]
&\approx&\int_{\vert x-(B_s-z)\vert\le a\}}
{1\over\vert x-(B_s-z)\vert^p}
\omega(z+dx)\\[-2pt]
&\ge& (\delta+r)^{-p} \biggl({\log t\over\log\log t} \biggr)^{p/ (2-p)}
\omega\biggl\{x; \vert z+x\vert\le\delta
\biggl({\log\log t\over\log t} \biggr)^{1/(2-p)} \biggr\}.\nonumber
\end{eqnarray}
Write
\[
X_z=\omega\biggl\{y; \vert z+y\vert\le\delta
\biggl({\log\log t\over\log t} \biggr)^{1/(2-p)} \biggr\}.
\]
In view of (\ref{theorem-17}),
\begin{eqnarray*}
&&\E_0\exp\biggl\{\theta\int_0^t\ol{V}(B_s)\,ds \biggr\}\\[-2pt]
&&\qquad \succeq\exp\biggl\{(r+\delta)^{-p}\theta t \biggl({\log t\over\log\log t}
\biggr)^{p/(2-p)}\max_z X_z\\[-2pt]
&&\hphantom{\qquad \succeq\exp\biggl\{}\hspace*{45pt}
{}-\lambda_d r^{-2}t \biggl({\log t\over\log\log t} \biggr)^{2/(2-p)} \biggr\}.
\end{eqnarray*}
Similarly to (\ref{theorem-21}),
\[
\lim_{t\to\infty}{\log\log t\over\log t}\max_z X_z={d(2-p)\over3-p}
\qquad \mbox{a.s.}
\]
Since $\delta>0$ can be arbitrarily small, the optimal pick
\[
r= \biggl({2\lambda_d(3-p)\over dp(2-p)\theta} \biggr)^{1/(2-p)}
\]
leads to the lower bound
\begin{eqnarray*}
&&\liminf_{t\to\infty}{1\over t} \biggl({\log\log t\over\log t} \biggr)^{2/(2-p)}
\E_0\exp\biggl\{\theta\int_0^t\ol{V}(B_s)\,ds \biggr\}\\[-2pt]
&&\qquad \ge{1\over2} \biggl({p\over2\lambda_d} \biggr)^{p/(2-p)}(2-p)^{(4-p)/(2-p)}
\biggl({d\theta\over3-p} \biggr)^{2/(2-p)}\qquad  \mbox{a.s.}
\end{eqnarray*}
This bound is sharp in rate in comparison with (\ref{theorem-11}).
Due to a lack of information on the value\vadjust{\goodbreak}
of $\sigma(d, p)$, we are not able to compare the constants on the
right-hand sides. However, it looks unlikely that the constant obtained here
would match the one in (\ref{theorem-11}). In addition, the argument
given in Sections \ref{u} and \ref{l}
shows that the constant $r>0$
should be arbitrarily large for the accuracy requested by (\ref{theorem-11}).

In summary, the simple strategies given above provide some heuristic pictures
on the behavior patterns of both Brownian particles and the Poisson
field and can be made rigorous for (\ref{intro-18}) and (\ref{intro-19}),
but fall short of the accuracy demanded by (\ref{theorem-10}) and
(\ref{theorem-11}). Some harder computation
on the tail estimates
for Poisson integrals is needed for the main theorems in this
paper.\looseness=-1

We now comment on the methods used in this paper.
The Feynman--Kac
formula is essential in this paper
for tracking the principal eigenvalues.
Among others, the
ingenious approach developed in \cite{G-K} and \cite{G-K-M},
which allows one to bound the principal eigenvalue over a
large domain by the maximal of the principal eigenvalues over
the sub-domains, plays a key role in our argument for the upper bound.
With this approach, we reduce the problem essentially to the tail estimate
of the random Dirichlet form
%
\begin{equation}\label{theorem-15}
\sup_{g\in\mathcal{F}_d(B(0,r\epsilon^{1/d}))} \!\biggl\{\!\pm\theta
\int_{B(0,r\epsilon^{1/d})}\!\ol{V}(x)g^2(x)\,dx\!-\!{1\over2}
\int_{B(0,r\epsilon^{1/d})}\!\vert\nabla g(x)\vert^2\,dx \!\biggr\},\hspace*{-38pt}
\end{equation}
where for any domain $D\subset\R^d$, $\mathcal{F}_d(D)$ is defined as
the set of the smooth functions $g$ on $D$ with $\|g\|_{\mathcal{L}^2(D)}=1$
and $g(\partial D)=0$, the constant $r>0$ is large but fixed and associated
with the critical radius $r_t$ posted in (\ref{theorem-18}), the parameter
$\epsilon$ is given as follows:
%
\begin{equation}\label{theorem-15'}\quad
\epsilon= \cases{
(\log t)^{-{(d-p)/2}},&\quad  in the proof of Theorem \ref{th-2},\vspace*{3pt}\cr
\displaystyle \biggl({\log\log t\over\log t} \biggr)^{d/(2-p)},&\quad
in the proof of Theorem \ref{th-3}.
}
\end{equation}

Another important idea adopted in this paper is the Poisson field rescaling.
In his proof of (\ref{intro-18}), Sznitman (\cite{Sznitman}, Chapter 4)
reduces the problem to the investigation of
the ``enlarged obstacles''
\[
\omega((\log t)^{1/d}\,dx ).
\]
It is worth pointing out that
the choice of the rescaling factor $(\log t)^{1/d}$ links to the critical
radius $r_t$ posted in (\ref{theorem-18}). What we
confront here are the ``contracted obstacles''
$\omega(\epsilon \,dx)-\epsilon\, dx$
with $\epsilon>0$ given
in (\ref{theorem-15'}).
Under the substitution
$g(x)\mapsto\epsilon^{-1/2}g(\epsilon^{-1/d}x)$
and by Fubini's theorem, the variation in (\ref{theorem-15}) is equal to
\begin{eqnarray*}
&&\sup_{g\in\mathcal{F}_d(B(0,r))} \biggl\{\pm\theta\epsilon^{-p/ d}
\int_{\R^d} \biggl[\int_{B(0,r)}{g^2(y)\over\vert y-x\vert^p}\,dy \biggr]
[\omega(\epsilon \,dx)-\epsilon \,dx ]\\[-3pt]
&&\hspace*{157pt}
{}-{\epsilon^{-2/d}\over2}
\int_{B(0,r)}\vert\nabla g(x)\vert^2\,dx \biggr\}.\vadjust{\goodbreak}
\end{eqnarray*}
The tail probabilities of the compensated Poisson integral appearing
here will be the main topic of the next section.

In comparison to the existing literature
such as \cite{BK1,G-K-M,HKM,HKS},
perhaps the most substantial difference comes from the fact that
in these works, the logarithmic moment generating
function (or the fractional logarithmic moment generating
function)
\[
H(\gamma)\equiv\log\E\exp\{\gamma V(0) \}
\]
exists. As a matter of fact, it is the logarithmic moment generating
function~$H(\gamma)$ (or the fractional logarithmic moment generating
function)
that plays a decisive role in these publications
in determining the asymptotics for
\begin{eqnarray*}
&&\log\E_0\exp\biggl\{\theta\int_0^tV(B_s)\,ds \biggr\} \quad\mbox{and}\\
&&\log\E\otimes\E_0\exp\biggl\{\theta\int_0^tV(B_s)\,ds \biggr\}\qquad  (t\to\infty)
\end{eqnarray*}
through some well-developed algorithms. Unfortunately,
this is not our case. Indeed, we have that
\[
\E\ol{V}^2(0)=\int_{\R^d}{dx\over\vert x\vert^{2p}}=\infty.
\]
Additional challenges we confront are the local unboundedness of $\ol{V}(x)$,
and the loss of monotonicity
of Poisson integrals due to renormalization.

The rest of the paper is organized as follows. In Section \ref{poisson}
we establish
the large deviations for a group of Poisson integrals with respect to
the contracted renormalized Poisson field. In
Section \ref{FK}, some explicit bounds for the Feynman--Kac
formula are established for later application. The upper bounds and
the lower bounds for our main theorems are proved in Sections~\ref{u}
and~\ref{l}, respectively. These bounds are established
simultaneously
for Theorems \ref{th-2} and  \ref{th-3}.
In Section \ref{A}, some identities
for the relevant integrals and variations are established.

\section{Large deviations for Poisson integrals}\label{poisson}

The functions
\[
\psi(\lambda)=e^{-\lambda}-1+\lambda\quad\mbox{and}\quad
\Psi(\lambda)=e^\lambda-1-\lambda\qquad (\lambda\ge0)
\]
appear frequently in this section. It is easy to see that
$\psi(\lambda)$ and $\Psi(\lambda)$ are nonnegative, increasing
and convex on $[0,\infty)$ with $\psi(0)=\Psi(0)=0$.
In addition, $\psi(\cdot)\le\Psi(\cdot)$ on $[0,\infty)$. According to
Lemma \ref{psi-1},
%
\begin{equation}\label{poisson-0}
\int_{\R^d}\psi\biggl({1\over\vert x\vert^p} \biggr)\,dx
=\omega_d{p\over d-p}\Gamma\biggl({2p-d\over p} \biggr)
\end{equation}
when $d/2<p<d$.\vadjust{\goodbreak}

The function $\Psi(\vert\cdot\vert^{-p})$ is not integrable on $\R^d$.
Under $p>d/2$, however,
\[
\int_{\{\vert x\vert\ge c\}}\Psi\biggl({1\over\vert x\vert^p}
\biggr)\,dx<\infty,
\qquad c>0.
\]

Throughout this section, $D\subset\R^d$ is a fixed bounded open set. Write
%
\begin{equation}\label{poisson-1}
\mathcal{G}_d(D)= \bigl\{g\in W^{1,2}(D); \|g\|^2_{\mathcal{L}^2(D)}+
\tfrac{1}{2}\|\nabla g\|^2_{\mathcal{L}^2(D)}=1 \bigr\},
\end{equation}
where $W^{1,2}(D)$ is the Sobolev space over $D$, defined to be the
closure of the inner product space
consisting of the infinitely differentiable functions compactly
supported in
$D$ under the Sobolev norm
\[
\|g\|_H= \bigl\{\|g\|^2_{\mathcal{L}^2(D)}+\|\nabla g\|^2_{\mathcal{L}^2(D)}
\bigr\}^{1/2}.
\]

To reserve
continuity we adopt a smooth truncation to the shape function.
Let the smooth function
$\alpha\dvtx\R^+\longrightarrow
[0,1]$ satisfy the following properties: $\alpha(\lambda)=1$ on $[0,1]$,
$\alpha(\lambda)=0$ for $\lambda\ge3$ and $-1\le\alpha'(\lambda)\le0$.

For $a>0$ and $\epsilon>0$, define
\begin{eqnarray*}
K^{(0)}_{a,\epsilon}(x)&=&{1\over\vert x\vert^p}
\alpha\bigl(a^{-1}\epsilon^{(2+d-p)/ (d(d-p))}\vert x\vert\bigr),
\\
K^{(1)}_{a,\epsilon}(x)
&=&{1\over\vert x\vert^p}\alpha(a^{-1}(\log\epsilon^{-1})^{-1/p}
\vert x\vert)
\end{eqnarray*}
and
\begin{eqnarray*}
L^{(0)}_{a,\epsilon}(x)&=&{1\over\vert x\vert^p} \bigl\{1-
\alpha\bigl(a^{-1}\epsilon^{(2+d-p)/(d(d-p))}
\vert x\vert\bigr) \bigr\},\\
L^{(1)}_{a,\epsilon}(x)
&=&{1\over\vert x\vert^p} \{1-\alpha(a^{-1}(\log\epsilon^{-1})^{-1/p}
\vert x\vert) \}
\end{eqnarray*}
and
%
\begin{eqnarray}\label{poisson-4}
G^{(i)}_{a,\epsilon}(g)=\int_{\R^d} \biggl[\int_DK^{(i)}_{a,\epsilon
}(y-x)g^2(y)\,dy \biggr]
[\omega(\epsilon\,dx)-\epsilon\,dx ],\nonumber\\[-8pt]\\[-8pt]
\eqntext{g\in\mathcal{G}_d(D),\
i=0,1,}
\\\label{poisson-5}
F^{(i)}_{a,\epsilon}(g)=\int_{\R^d} \biggl[\int_DL^{(i)}_{a,\epsilon
}(y-x)g^2(y)\,dy \biggr]
[\omega(\epsilon\,dx)-\epsilon\,dx ],\nonumber\\[-8pt]\\[-8pt]
\eqntext{g\in\mathcal{G}_d(D),\ i=0,1.}
\end{eqnarray}

Write
%
\begin{equation}\label{poisson-2}
\zeta_\epsilon(g)=\int_{\R^d} \biggl[\int_D{g^2(y)\over\vert y-x\vert^p}\,dy \biggr]
[\omega(\epsilon\,dx)-\epsilon\,dx ],\qquad  g\in\mathcal{G}_d(D).
\end{equation}

The main theorems in this section are the large deviations for the Poisson
integrals indexed by $\mathcal{G}_d(D)$.

\begin{theorem}\label{poisson-6} Assume that $d/2<p<d$.
For any $a>0$ and $\gamma>0$,
%
\begin{eqnarray} \label{poisson-7}
\lim_{a\to\infty}\limsup_{\epsilon\to0^+}\epsilon^{2/(d-p)}\log
\P\Bigl\{\sup_{g\in\mathcal{G}_d(D)}\bigl\vert F^{(0)}_{a,\epsilon}(g)\bigr\vert
\ge\gamma\epsilon^{-{(2-p)/d}}\Bigr \}&=&-\infty,\hspace*{-10pt}
\\[-2pt]
\label{poisson-8}
\qquad\quad  \liminf_{a\to\infty}\liminf_{\epsilon\to0^+}\epsilon^{2/(d-p)}\log
\P\Bigl\{\inf_{g\in\mathcal{G}_d(D)}G^{(0)}_{a,\epsilon}(g)
\le-\gamma\epsilon^{-{(2-p)/d}} \Bigr\}&\ge&-I_D(\gamma),\hspace*{-10pt}
\\[-2pt]
\label{poisson-9}
\lim_{\epsilon\to0^+}\epsilon^{2/(d-p)}\log
\P\Bigl\{\inf_{g\in\mathcal{G}_d(D)}\zeta_\epsilon(g)
\le-\gamma\epsilon^{-{(2-p)/d}} \Bigr\}&=&-I_D(\gamma),\hspace*{-10pt}
\end{eqnarray}
where
\begin{eqnarray} \label{poisson-10}
I_D(\gamma)&=& \biggl({\gamma(d-p)\over d} \biggr)^{d/(d-p)}
\biggl(\omega_d
\Gamma\biggl({2p-d\over p} \biggr)
\biggr)^{-{p/(d-p)}}\nonumber\\[-9pt]\\[-9pt]
&&{}\times\Bigl(\sup_{g\in\mathcal{G}_d(D)}\|g\|_{\mathcal{L}^2(D)} \Bigr)^{-{2d/(d-p)}}.\nonumber
\end{eqnarray}
\end{theorem}

Write $ l(\epsilon)=\epsilon^{-{(2-p)/d}}\log{1\over\epsilon}$,
and
%
\begin{equation}\label{poisson-13}
\rho^*_D=\sup_{g\in\mathcal{G}_d(D)}\sup_{x\in\R^d}
\int_D{g^2(y)\over\vert y-x\vert^p}\,dy.
\end{equation}
The finiteness of $\rho^*_D$ can be seen from
(\ref{theorem-6'}) in Lemma \ref{th-1} and from (\ref{poisson-19}).

\begin{theorem}\label{poisson-14} Assume $d/2<p<\min\{2, d\}$.
For any $a>0$ and $\gamma>0$,
%
\begin{eqnarray} \label{poisson-15}
\qquad \lim_{\epsilon\to0^+}{1\over l(\epsilon)}
\log\P\Bigl\{\sup_{g\in\mathcal{G}_d(D)}\bigl\vert F^{(1)}_{a,\epsilon
}(g)\bigr\vert\ge
\gamma\epsilon^{-{(2-p)/d}} \Bigr\}&=&-\infty,
\\[-2pt]
\label{poisson-16}
\lim_{\epsilon\to0^+}{1\over l(\epsilon)}
\log\P\Bigl\{\sup_{g\in\mathcal{G}_d(D)}G^{(1)}_{a,\epsilon}(g)
\ge\gamma\epsilon^{-{(2-p)/d}} \Bigr\}&=& -{2+d-p\over d\rho^*_D}\gamma,
\\[-2pt]
\label{poisson-17}
\lim_{\epsilon\to0^+}{1\over l(\epsilon)}
\log\P\Bigl\{\sup_{g\in\mathcal{G}_d(D)}\zeta_\epsilon(g)
\ge\gamma\epsilon^{-{(2-p)/d}} \Bigr\}&=&-{2+d-p\over d\rho^*_D}
\gamma.
\end{eqnarray}
\end{theorem}

Write
\[
\ol{V}^{(i)}_{a,\epsilon}(x)=\int_{\R^d}L_{a,\epsilon}^{(i)}(y-x)
[\omega(\epsilon \,dy)-\epsilon \,dy ],\qquad
x\in D, i=0, 1.
\]

Our approach relies on the following lemma.

\begin{lemma}\label{poisson-19'} For any $a>0$ and $\theta>0$,
%
\begin{eqnarray}\label{poisson-22'}
&&\lim_{\epsilon\to0^+}\epsilon^{2/(d-p)}\log
\E\exp\Bigl\{-\theta\epsilon^{-{p(2+d-p)/(d(d-p))}}
\inf_{x\in D}\ol{V}^{(0)}_{a,\epsilon}(x) \Bigr\}\nonumber\\[-9pt]\\[-9pt]
&&\qquad =\int_{\R^d}
\psi\biggl(\theta{1-\alpha(a^{-1}\vert x\vert)\over\vert x\vert^p} \biggr)\,dx,\nonumber
\\
\label{poisson-22''}
&&\lim_{\epsilon\to0^+}\epsilon^{2/(d-p)}\log
\E\exp\Bigl\{\theta\epsilon^{-{p(2+d-p)/(d(d-p))}}
\sup_{x\in
D}\bigl\vert\ol{V}^{(0)}_{a,\epsilon}(x)\bigr\vert\Bigr\}\nonumber\\[-8pt]\\[-8pt]
&&\qquad =\int_{\R^d}
\Psi\biggl(\theta{1-\alpha(a^{-1}\vert x\vert)\over\vert x\vert^p} \biggr)\,dx,\nonumber
\\
\label{poisson-23}
&&\lim_{\epsilon\to0^+}{1\over l(\epsilon)}
\log\E\exp\biggl\{\theta \biggl(\log{1\over\epsilon} \biggr)
\sup_{x\in D}\bigl\vert\ol{V}^{(1)}_{a,\epsilon}(x)\bigr\vert\biggr\}=0.
\end{eqnarray}
\end{lemma}

\begin{pf}
Notice that for any $a>0$, $\theta>0$ and $x\in D$,
\begin{eqnarray*}
&&\E\exp\bigl\{-\theta\epsilon^{-{p(2+d-p)/(d(d-p))}}
\ol{V}^{(0)}_{a,\epsilon}(x) \bigr\}\\
&&\qquad =\E\exp\bigl\{-\theta
\epsilon^{-{p(2+d-p)/(d(d-p))}}\ol{V}^{(0)}_{a,\epsilon}(0) \bigr\}\\
&&\qquad =\exp\biggl\{\epsilon\int_{\R^d}\psi\bigl(\theta\epsilon^{-{p(2+d-p)/(d(d-p))}}
L^{(0)}_{a,\epsilon}(x) \bigr)\,dx \biggr\}\\
&&\qquad =\exp\biggl\{\epsilon^{2/(d-p)}\int_{\R^d}
\psi\biggl(\theta{1-\alpha(a^{-1}\vert x\vert)\over\vert x\vert^p} \biggr)\,dx \biggr\}.
\end{eqnarray*}
Similarly,
\begin{eqnarray*}
&&\E\exp\bigl\{\theta\epsilon^{-{p(2+d-p)/(d(d-p))}}
\ol{V}^{(0)}_{a,\epsilon}(x) \bigr\}\\
&&\qquad =\exp\biggl\{\epsilon^{2/(d-p)}\int_{\R^d}
\Psi\biggl(\theta{1-\alpha(a^{-1}\vert x\vert)\over\vert x\vert^p} \biggr)\,dx \biggr\}.
\end{eqnarray*}
In view of the fact that $\psi(\cdot)\le\Psi(\cdot)$ on $[0,\infty)$,
we conclude that
\begin{eqnarray}\label{lemma-1}
&&\lim_{\epsilon\to0^+}\epsilon^{2/(d-p)}
\log\E\exp\bigl\{-\theta\epsilon^{-{p(2+d-p)/(d(d-p))}}
\ol{V}^{(0)}_{a,\epsilon}(x) \bigr\}\nonumber\\[-8pt]\\[-8pt]
&&\qquad =\int_{\R^d}
\psi\biggl(\theta{1-\alpha(a^{-1}\vert x\vert)\over\vert x\vert^p} \biggr)\,dx,\nonumber
\\
\label{lemma-2}
&&\lim_{\epsilon\to0^+}\epsilon^{2/(d-p)}
\log\E\exp\bigl\{\theta\epsilon^{-{p(2+d-p)/(d(d-p))}}
\bigl\vert\ol{V}^{(0)}_{a,\epsilon}(x)\bigr\vert\bigr\}\nonumber\\[-8pt]\\[-8pt]
&&\qquad =\int_{\R^d}
\Psi\biggl(\theta{1-\alpha(a^{-1}\vert x\vert)\over\vert x\vert^p} \biggr)\,dx.\nonumber
\end{eqnarray}

A similar computation also leads to
%
\begin{equation}\label{lemma-3}
\lim_{\epsilon\to0^+}{1\over l(\epsilon)}
\log\E\exp\biggl\{\theta \biggl(\log{1\over\epsilon} \biggr)
\bigl\vert\ol{V}^{(1)}_{a,\epsilon}(x)\bigr\vert\biggr\}=0,\qquad x\in D.
\end{equation}

All we need is to take supremum over $x\in D$ in the exponent
on the left-hand
sides of (\ref{lemma-1}), (\ref{lemma-2}) and (\ref{lemma-3})
and push the supremum through the expectation. Due to similarity,
we only carry out this algorithm to (\ref{lemma-1}) and
(\ref{lemma-2}).\vadjust{\goodbreak}

By the boundedness of $D$, we may assume that $D=(-b, b)^d$ for some $b>0$.
Let $h>0$ be a constant which will be later specified, and let
\[
\gamma={2+d-p\over d-p}+h.
\]
By integration substitution
\begin{eqnarray}\label{lemma-4}
\ol{V}^{(0)}_{a,\epsilon}(x)&=&\epsilon^{-ph/d}\int_{\R^d}
\widetilde{L}^{(0)}_{a,\epsilon} (y-\epsilon^{-h/d}x )
[\omega(\epsilon^{1+h}\,dy)-\epsilon^{1+h}\,dy
]\nonumber\\[-9pt]\\[-9pt]
&=&\epsilon^{-ph/d}H_\epsilon(\epsilon^{-h/d}x),\nonumber
\end{eqnarray}
where
\[
\widetilde{L}^{(0)}_{a,\epsilon}(x)={1\over\vert x\vert^p}
\{1-\alpha(a^{-1}\epsilon^{\gamma/d}\vert x\vert) \}
\]
and
\[
H_\epsilon(x)=\int_{\R^d}
\widetilde{L}^{(0)}_{a,\epsilon} (z-x )
[\omega(\epsilon^{1+h}\,dz)-\epsilon^{1+h}\,dz ].
\]

For any $x, y\in D$ with $x\not=y$, and $\theta>0$,
\begin{eqnarray*}
&&\E\exp\biggl\{\theta\epsilon^{-p\gamma/d}{H_\epsilon(x)-H_\epsilon(y)\over
\vert x-y\vert} \biggr\}\\[-2pt]
&&\qquad =\exp\biggl\{\epsilon^{1+h}\int_{\R^d}
\Psi\biggl(\epsilon^{-p\gamma/d}{\theta\over\vert x-y\vert}
\bigl(\widetilde{L}^{(0)}_{a,\epsilon}(z-x)-\widetilde{L}^{(0)}_{a,\epsilon}(z-y)
\bigr) \biggr)\,dz \biggr\}.
\end{eqnarray*}
Switching $x$ and $y$, one has
\begin{eqnarray*}
&&\E\exp\biggl\{\theta\epsilon^{-p\gamma/d}{
\vert H_\epsilon(x)-H_\epsilon(y)\vert\over
\vert x-y\vert} \biggr\}\\[-2pt]
&&\qquad \le2\exp\biggl\{\epsilon^{1+h}\int_{\R^d}
\Psi\biggl(\epsilon^{-p\gamma/d}{\theta\over\vert x-y\vert}
\bigl\vert\widetilde{L}^{(0)}_{a,\epsilon}(z-x)
-\widetilde{L}^{(0)}_{a,\epsilon}(z-y)
\bigr\vert\biggr)\,dz \biggr\}.
\end{eqnarray*}
By integration substitution,
\begin{eqnarray*}
&&\int_{\R^d}
\Psi\biggl(\epsilon^{-p\gamma/d}{\theta\over\vert x-y\vert}
\bigl\vert\widetilde{L}^{(0)}_{a,\epsilon}(z-x)-\widetilde
{L}^{(0)}_{a,\epsilon}(z-y)
\bigr\vert\biggr)\,dz\\[-2pt]
&&\qquad =\epsilon^{-\gamma}\int_{\R^d}\Psi\biggl({\theta\over\vert x-y\vert}
\vert L_a (z-\epsilon^{\gamma/d}x )-
L_a (z-\epsilon^{\gamma/d}y ) \vert\biggr)\,dz,
\end{eqnarray*}
where
\[
L_a(z)={1-\alpha(a^{-1}\vert z\vert)\over\vert z\vert^p}.
\]
By the mean value theorem, there is a $C_a>0$ such that when $\epsilon>0$
is sufficiently small,
\[
\vert L_a (z-\epsilon^{\gamma/d}x )-
L_a (z-\epsilon^{\gamma/d}y ) \vert\le C_a{\epsilon^{\gamma/d}
\vert x-y\vert\over\vert z\vert^p}
1_{\{\vert z\vert\ge C_a^{-1}\}}, \qquad x, y\in D.\vadjust{\goodbreak}
\]

Summarizing what we have,
\begin{eqnarray*}
&&\E\exp\biggl\{\theta\epsilon^{-p\gamma/d}{
\vert H_\epsilon(x)-H_\epsilon(y)\vert\over
\vert x-y\vert} \biggr\}\\
&&\qquad \le2\exp\biggl\{\epsilon^{-{2/(d-p)}}
\int_{\{\vert x\vert\ge C_a^{-1}\}}\Psi\biggl({C_a\theta\epsilon^{\gamma/d}
\over\vert z\vert^p} \biggr)\,dz \biggr\}\\
&&\qquad =2\exp\biggl\{\epsilon^{{\gamma/ p}-{2/ (d-p)}}
\int_{\{\vert x\vert\ge C_a^{-1}\epsilon^{-\gamma/(dp)}\}}
\Psi\biggl({C_a\theta
\over\vert z\vert^p} \biggr)\,dz \biggr\}.
\end{eqnarray*}
Let $h>0$ satisfy that
\[
h\ge{3p-d-2\over d-p}\quad  \mbox{or}\quad
{\gamma\over p}-{2\over d-p}\ge0.
\]
Then for any $\theta>0$ the quantity
\[
\mathop{\sup_{x,y\in D}}_{x\not=y}
\E\exp\biggl\{\theta\epsilon^{-p\gamma/d}{
\vert H_\epsilon(x)-H_\epsilon(y)\vert\over
\vert x-y\vert} \biggr\}
\]
is bounded uniformly for small $\epsilon>0$. Thus (\cite{Chen}, Theorem D-6), for any $\theta>0$,
%
\begin{equation}\label{lemma-5}
\lim_{\delta\to0^+}\limsup_{\epsilon\to0^+}
\E\exp\Bigl\{\theta\epsilon^{-p\gamma/d}
\sup_{\vert x-y\vert\le\delta}\vert H_\epsilon(x)-H_\epsilon(y)\vert\Bigr\}=1.
\end{equation}

On the other hand, for any $x\in\R^d$ and $\theta>0$,
\begin{eqnarray*}
\E\exp\{\pm\theta\epsilon^{-p\gamma/d}H_\epsilon(x) \}
&=&\E\exp\{\pm\theta\epsilon^{-p\gamma/d}H_\epsilon(0) \}\\
&=&\E\exp\bigl\{\pm\theta\epsilon^{-{p(2+d-p)/(d(d-p))}}
\ol{V}^{(0)}(0) \bigr\},
\end{eqnarray*}
where the last step follows from (\ref{lemma-4}). By (\ref{lemma-1}) and
(\ref{lemma-2}), therefore, for any $x\in D$
\begin{eqnarray*}
\lim_{\epsilon\to0^+}\epsilon^{2/(d-p)}
\log\E\exp\{-\theta\epsilon^{-p\gamma/d}H_\epsilon(x) \}
&=&\int_{\R^d}
\psi\biggl(\theta{1-\alpha(a^{-1}\vert y\vert)\over\vert y\vert^p} \biggr)\,dy,
\\
\lim_{\epsilon\to0^+}\epsilon^{2/(d-p)}
\log\E\exp\{\theta\epsilon^{-p\gamma/d}\vert H_\epsilon(x)\vert\}
&=&\int_{\R^d}
\Psi\biggl(\theta{1-\alpha(a^{-1}\vert y\vert)\over\vert y\vert^p} \biggr)\,dy.
\end{eqnarray*}
Combine them with (\ref{lemma-5}). A standard argument of
exponential approximation leads to
\begin{eqnarray}\label{lemma-6}
&&\lim_{\epsilon\to0^+}\epsilon^{2/ (d-p)}
\log\E\exp\Bigl\{-\theta\epsilon^{-p\gamma/d}\inf_{x\in D}H_\epsilon(x)
\Bigr\}\nonumber\\[-8pt]\\[-8pt]
&&\qquad =\int_{\R^d}
\psi\biggl(\theta{1-\alpha(a^{-1}\vert x\vert)\over\vert x\vert^p} \biggr)\,dx,\nonumber
\\
\label{lemma-7}
&&\lim_{\epsilon\to0^+}\epsilon^{2/ (d-p)}
\log\E\exp\Bigl\{\theta\epsilon^{-p\gamma/d}\sup_{x\in D}
\vert H_\epsilon(x)\vert\Bigr\}\nonumber\\[-8pt]\\[-8pt]
&&\qquad =\int_{\R^d}
\Psi\biggl(\theta{1-\alpha(a^{-1}\vert x\vert)\over\vert x\vert^p}
\biggr)\,dx.\nonumber
\end{eqnarray}

Recall that $D=(-b, b)^d$. Using (\ref{lemma-4}),
\begin{eqnarray*}
-\inf_{x\in D}V_{a,\epsilon}^{(0)}(x)&=&-\epsilon^{-ph/d}\inf_{x\in
\epsilon^{-h/d}D}
H_\epsilon(x)\\
&\le&\epsilon^{-ph/d}
\max_{z\in b\Z^d\cap\epsilon^{-h/d}D} \Bigl\{-\inf_{x\in z+D}H_\epsilon(x)\Bigr\}.
\end{eqnarray*}
By the fact that the random variables
\[
\inf_{x\in z+D}H_\epsilon(x);\qquad  z\in b\Z^d\cap\epsilon^{-h/d}D,
\]
are identically distributed,
\begin{eqnarray*}
&&\E\exp\Bigl\{-\theta\epsilon^{-{p(2+d-p)/(d(d-p))}}
\inf_{x\in D}\ol{V}^{(0)}_{a,\epsilon}(x) \Bigr\}\\
&&\qquad \le\# \{ b\Z^d\cap\epsilon^{-h/d}D \}
\E\exp\Bigl\{-\theta\epsilon^{-p\gamma/d}\inf_{x\in D}H_\epsilon(x) \Bigr\}.
\end{eqnarray*}
Consequently from (\ref{lemma-6}),
\begin{eqnarray*}
&&\limsup_{\epsilon\to0^+}\epsilon^{2/(d-p)}
\log\E\exp\Bigl\{-\theta\epsilon^{-{p(2+d-p)/(d(d-p))}}
\inf_{x\in D}\ol{V}^{(0)}_{a,\epsilon}(x) \Bigr\}\\
&&\qquad \le\int_{\R^d}
\Psi\biggl(\theta{1-\alpha(a^{-1}\vert x\vert)\over\vert x\vert^p} \biggr)\,dx.
\end{eqnarray*}
In view of (\ref{lemma-1}), we have proved (\ref{poisson-22'}).

Assertion (\ref{poisson-22''}) follows from (\ref{lemma-2}) and
(\ref{lemma-7}) in the same way.
\end{pf}

\subsection{\texorpdfstring{Proof of Theorem \protect\ref{poisson-6}}{Proof of Theorem 3.1}}

Let $\theta>0$ be fixed but arbitrary. By (\ref{poisson-22''}) and
the inequality
\begin{eqnarray*}
\sup_{g\in\mathcal{G}_d(D)}\vert F_{a,\epsilon}(g)\vert&\le&\sup_{g\in
\mathcal{G}_d(D)}
\int_D\bigl\vert\ol{V}_{a,\epsilon}^{(0)}(x)\bigr\vert g^2(x)\,dx\\
&\le& \Bigl(\sup_{g\in\mathcal{G}_d(D)}\|g\|_{\mathcal{L}^2(D)} \Bigr)^2
\sup_{x\in D}\bigl\vert\ol{V}_{a,\epsilon}^{(0)}(x)\bigr\vert,
\end{eqnarray*}
we have
\begin{eqnarray*}
&&\limsup_{\epsilon\to0^+}\epsilon^{2/ (d-p)}\log
\E\exp\Bigl\{\theta\epsilon^{-{p(2+d-p)/ (d(d-p))}}
\sup_{g\in\mathcal{G}_d(D)}\bigl\vert F^{(0)}_{a,\epsilon}(g)\bigr\vert\Bigr\}\\
&&\qquad\le\int_{\R^d}\Psi\biggl( \Bigl(\sup_{g\in\mathcal{G}_d(D)}
\|g\|_{\mathcal{L}^2(D)} \Bigr)^2\theta
{1-\alpha(a^{-1}\vert x\vert)
\over\vert x\vert^p} \biggr)\,dx.
\end{eqnarray*}
Consequently,
\[
\lim_{a\to\infty}\limsup_{\epsilon\to0^+}\epsilon^{2/ (d-p)}\log
\E\exp\Bigl\{\theta\epsilon^{-{p(2+d-p)/(d(d-p))}}
\sup_{g\in\mathcal{G}_d(D)}\bigl\vert F^{(0)}_{a,\epsilon}(g)\bigr\vert\Bigr\}
=0.
\]
Therefore, (\ref{poisson-7}) follows from a standard application
of Chebyshev's inequality.

We now prove (\ref{poisson-9}). For any $g\in\mathcal{G}_d(D)$,
\begin{eqnarray*}
&&\E\exp\bigl\{-\theta\epsilon^{-{p(2+d-p)/(d(d-p))}}
\zeta_\epsilon(g) \bigr\}\\[-2pt]
&&\qquad =\exp\biggl\{\epsilon\int_{\R^d}\psi\biggl(\theta
\epsilon^{-{p(2+d-p)/(d(d-p))}}\int_D{1\over\vert
y-x\vert^p}g^2(y)\,dy \biggr)\,dx \biggr\}.
\end{eqnarray*}
Given $\delta>0$,
%
\begin{eqnarray}
\hspace*{-3pt}&&\int_{\R^d}\psi\biggl(\theta
\epsilon^{-{p(2+d-p)/(d(d-p))}}\int_D{1\over\vert
y-x\vert^p}g^2(y)\,dy \biggr)\,dx\nonumber \\[-2pt]
\hspace*{-3pt}&&\qquad \ge\int_{\{\vert x\vert\ge\delta\epsilon^{-{(2+d-p)/(d(d-p))}}\}}
\psi\biggl(\theta\epsilon^{-{p(2+d-p)/(d(d-p))}}\int_D{1\over\vert
y-x\vert^p}g^2(y)\,dy \biggr)\,dx\nonumber \\[-2pt]
\hspace*{-3pt}&&\qquad \ge\int_{\{\vert x\vert\ge\delta\epsilon^{-{(2+d-p)/(d(d-p))}}\}}
\psi\biggl(\theta
\epsilon^{-{p(2+d-p)/(d(d-p))}}{ (1+o(1) )\over\vert x\vert^p}
\|g\|_{\mathcal{L}^2(D)}^2 \biggr)\,dx\nonumber \\[-2pt]
\hspace*{-3pt}&&\qquad = \bigl(1+o(1) \bigr)\theta^{d/p}\|g\|_{\mathcal{L}^2(D)}^{2d/p}
\epsilon^{-{(2+d-p)/(d-p)}}\int_{\{\vert x\vert\ge(1+o(1))\delta\}}
\psi\biggl({1\over\vert x\vert^p} \biggr)\,dx\nonumber \\[-2pt]
\hspace*{-3pt}\eqntext{(\epsilon\to0^+).}
\end{eqnarray}

Since $\delta$ can be arbitrarily
small, we have
%
\begin{eqnarray}\label{poisson-25}
&&\liminf_{\epsilon\to0^+}\epsilon^{2/ (d-p)}\log
\E\exp\Bigl\{-\theta\epsilon^{-{p(2+d-p)/ (d(d-p))}}\inf_{g\in\mathcal{G}_d(D)}
\zeta_\epsilon(g) \Bigr\}\nonumber\\[-2pt]
&&\qquad\ge\theta^{d/p} \Bigl(\sup_{g\in\mathcal{G}_d(D)}\|g\|_{\mathcal
{L}^2(D)} \Bigr)^{2d/p}
\int_{\R^d}
\psi\biggl({1\over\vert x\vert^p} \biggr)\,dx\\[-2pt]
&&\qquad=\theta^{d/p} \Bigl(\sup_{g\in\mathcal{G}_d(D)}\|g\|_{\mathcal
{L}^2(D)} \Bigr)^{2d/p}
{\omega_dp\over d-p}\Gamma\biggl({2p-d\over p} \biggr),\nonumber
\end{eqnarray}
where the last step follows from (\ref{poisson-0}).

On the other hand, for any $g\in\mathcal{G}_d(D)$ and $a>0$,
%
\begin{equation}\label{poisson-26}
\zeta_\epsilon(g)=G^{(0)}_{a,\epsilon}(g)+F^{(0)}_{a,\epsilon}(g).
\end{equation}
Notice that
\begin{eqnarray*}
G^{(0)}_{a,\epsilon}(g)&\ge&-\epsilon\int_{\R^d} \biggl[\int_D
K^{(0)}_{a,\epsilon}(y-x)g^2(y)\,dy \biggr]\,dx\\[-2pt]
&=&-\epsilon\|g\|_{\mathcal{L}^2(D)}^2\int_{\R^d}K^{(0)}_{a,\epsilon}(x)\,dx\\[-2pt]
&\ge&- Ca^{d-p}\epsilon^{-{(2-p)/d}}.
\end{eqnarray*}
Consequently,
%
\begin{eqnarray}\label{poisson-27}
 &&\limsup_{\epsilon\to0^+}\epsilon^{2/ (d-p)}\log
\E\exp\Bigl\{-\theta\epsilon^{-{p(2+d-p)/ (d(d-p))}}\inf_{g\in\mathcal{G}_d(D)}
\zeta_\epsilon(g) \Bigr\}\nonumber\hspace*{-35pt}\\
&&\quad\le C\theta a^{d-p}\hspace*{-35pt}\\
&&\qquad  {}+\limsup_{\epsilon\to0^+}\epsilon^{2/
(d-p)}\log
\E\exp\Bigl\{-\theta\epsilon^{-{p(2+d-p)/ (d(d-p))}}\inf_{g\in\mathcal{G}_d(D)}
F_{a,\epsilon}^{(0)} (g) \Bigr\}.\nonumber\hspace*{-35pt}
\end{eqnarray}

To deal with the right-hand side, notice that
\[
F^{(0)}_{a,\epsilon}(g)=
\int_D\ol{V}^{(0)}_{a,\epsilon}(x)g^2(x)\,dx
\ge\|g\|_{\mathcal{L}^2(D)}^2\inf_{x\in D}\ol{V}^{(0)}_{a,\epsilon}(x).
\]
Hence,
\[
\inf_{g\in\mathcal{G}_d(D)}F^{(0)}_{a,\epsilon}(g)=
\inf_{g\in\mathcal{G}_d(D)}\int_D\ol{V}^{(0)}_{a,\epsilon}(x)g^2(x)\,dx
\ge \Bigl(\sup_{g\in\mathcal{G}_d(D)}\|g\|_{\mathcal{L}^2(D)}^2 \Bigr)\inf_{x\in D}
\ol{V}^{(0)}_{a,\epsilon}(x)
\]
when $\inf_{x\in D}\ol{V}^{(0)}_{a,\epsilon}(x)\le0$, and
\[
\inf_{g\in\mathcal{G}_d(D)}F^{(0)}_{a,\epsilon}(g)=
\inf_{g\in\mathcal{G}_d(D)}\int_D\ol{V}^{(0)}_{a,\epsilon}(x)g^2(x)\,dx
\ge \Bigl(\inf_{g\in\mathcal{G}_d(D)}\|g\|_{\mathcal{L}^2(D)}^2 \Bigr)\inf_{x\in D}
\ol{V}^{(0)}_{a,\epsilon}(x)
\]
when $\inf_{x\in D}\ol{V}^{(0)}_{a,\epsilon}(x)> 0$. Thus,
\begin{eqnarray*}
&&\E\exp\Bigl\{-\theta\epsilon^{-{p(2+d-p)/ (d(d-p))}}\inf_{g\in\mathcal{G}_d(D)}
F_{a,\epsilon}^{(0)} (g) \Bigr\}\\
&&\qquad\le
\E\Bigl[\exp\Bigl\{-\theta\epsilon^{-{p(2+d-p)/ (d(d-p))}}
\sup_{g\in\mathcal{G}_d(D)}\|g\|_{\mathcal{L}^2(D)}^2\inf_{x\in D}
\ol{V}^{(0)}_{a,\epsilon}(x) \Bigr\};\\
&&\hspace*{238pt}
\inf_{x\in D}\ol{V}^{(0)}_{a,\epsilon}(x)\le0 \Bigr]\\
&&\qquad \quad {}+\E\Bigl[\exp\Bigl\{-\theta\epsilon^{-{p(2+d-p)/ (d(d-p))}}
\inf_{g\in\mathcal{G}_d(D)}\|g\|_{\mathcal{L}^2(D)}^2\inf_{x\in D}
\ol{V}^{(0)}_{a,\epsilon}(x) \Bigr\};\\
&&\hspace*{252pt}
\inf_{x\in D}\ol{V}^{(0)}_{a,\epsilon}(x)> 0 \Bigr]\\
&&\qquad\le
\E\exp\Bigl\{-\theta\epsilon^{-{p(2+d-p)/ (d(d-p))}}
\sup_{g\in\mathcal{G}_d(D)}\|g\|_{\mathcal{L}^2(D)}^2\inf_{x\in D}
\ol{V}^{(0)}_{a,\epsilon}(x) \Bigr\}\\
&&\qquad\quad{}+\E\exp\Bigl\{-\theta\epsilon^{-{p(2+d-p)/ (d(d-p))}}
\inf_{g\in\mathcal{G}_d(D)}\|g\|_{\mathcal{L}^2(D)}^2\inf_{x\in D}
\ol{V}^{(0)}_{a,\epsilon}(x) \Bigr\}.
\end{eqnarray*}
By (\ref{poisson-22'}) with $\theta$ being replaced by
\[
\theta\sup_{g\in\mathcal{G}_d(D)}\|g\|^2_{\mathcal{L}^2(D)}\quad
\mbox{and}\quad
\theta
\inf_{g\in\mathcal{G}_d(D)}\|g\|^2_{\mathcal{L}^2(D)},
\]
respectively,
\begin{eqnarray*}
&&\limsup_{\epsilon\to0^+}\epsilon^{2/(d-p)}\log
\E\exp\Bigl\{-\theta\epsilon^{-{p(2+d-p)/(d(d-p))}}\inf_{g\in\mathcal{G}_d(D)}
F^{(0)}_{a,\epsilon}(g) \Bigr\}\\
&&\qquad\le\int_{\R^d}
\psi\biggl({\theta (1-\alpha(a^{-1}\vert x\vert) )
\over\vert x\vert^p}\sup_{g\in\mathcal{G}_d(D)}\|g\|^2_{\mathcal
{L}^2(D)} \biggr)\,dx\\
&&\qquad\le\int_{\R^d}
\psi\biggl({\theta
\over\vert x\vert^p}\sup_{g\in\mathcal{G}_d(D)}\|g\|^2_{\mathcal
{L}^2(D)} \biggr)\,dx\\
&&\qquad=\theta^{d/p} \Bigl(\sup_{g\in\mathcal{G}_d(D)}\|g\|^2_{\mathcal
{L}^2(D)} \Bigr)^{2d/p}
\int_{\R^d}\psi\biggl({1\over\vert x\vert^p} \biggr)\,dx\\
&&\qquad=\theta^{d/p} \Bigl(\sup_{g\in\mathcal{G}_d(D)}\|g\|_{\mathcal
{L}^2(D)} \Bigr)^{2d/p}
{\omega_dp\over d-p}\Gamma\biggl({2p-d\over p} \biggr).
\end{eqnarray*}
Bringing this to (\ref{poisson-27}),
\begin{eqnarray*}
&&\limsup_{\epsilon\to0^+}\epsilon^{2/(d-p)}\log
\E\exp\Bigl\{-\theta\epsilon^{-{p(2+d-p)/(d(d-p))}}\inf_{g\in\mathcal{G}_d(D)}
\zeta_\epsilon(g) \Bigr\}\\
&&\qquad\le C\theta a^{d-p}+\theta^{d/p}
\Bigl(\sup_{g\in\mathcal{G}_d(D)}\|g\|_{\mathcal{L}^2(D)} \Bigr)^{2d/p}
{\omega_dp\over d-p}\Gamma\biggl({2p-d\over p}
\biggr).
\end{eqnarray*}
Letting $a\to0^+$ on the right-hand side leads to
\begin{eqnarray}\label{poisson-28}
&&\limsup_{\epsilon\to0^+}\epsilon^{2/ (d-p)}\log
\E\exp\Bigl\{-\theta\epsilon^{-{p(2+d-p)/ (d(d-p))}}\inf_{g\in\mathcal{G}_d(D)}
\zeta_\epsilon(g) \Bigr\}\nonumber\\[-8pt]\\[-8pt]
&&\qquad\le\theta^{d/p}
\Bigl(\sup_{g\in\mathcal{G}_d(D)}\|g\|_{\mathcal{L}^2(D)} \Bigr)^{2d/p}
{\omega_dp\over d-p}\Gamma\biggl({2p-d\over p} \biggr).\nonumber
\end{eqnarray}

The combination of (\ref{poisson-25}) and (\ref{poisson-28}) implies
(\cite{Chen}, Theorem 1.2.4) that
\begin{eqnarray*}
&&\lim_{\epsilon\to0^+}\epsilon^{2/(d-p)}\log
\P\Bigl\{\inf_{g\in\mathcal{G}_d(D)}\zeta_\epsilon(g)
\le-\gamma\epsilon^{-{(2-p)/d}} \Bigr\}\\
&&\qquad=-\sup_{\theta>0} \biggl\{\gamma\theta-\theta^{d/p}
\Bigl(\sup_{g\in\mathcal{G}_d(D)}\|g\|_{\mathcal{L}^2(D)} \Bigr)^{2d/p}
{\omega_dp\over d-p}\Gamma\biggl({2p-d\over p} \biggr) \biggr\}\\
&&\qquad =-I_D(\gamma).
\end{eqnarray*}

It remains to prove (\ref{poisson-8}). By (\ref{poisson-26}), for any
$\delta>0$,
\begin{eqnarray*}
&&\P\Bigl\{\inf_{g\in\mathcal{G}_d(D)}\zeta_\epsilon(g)
\le-(\gamma+\delta)\epsilon^{-{(2-p)/ d}} \Bigr\}\\
&&\qquad\le\P\Bigl\{\inf_{g\in\mathcal{G}_d(D)}G^{(0)}_{a,\epsilon}(g)
\le-\gamma\epsilon^{-{(2-p)/ d}} \Bigr\}\\
&&\qquad \quad {}+\P\Bigl\{\sup_{g\in\mathcal{G}_d(D)}\bigl\vert F^{(0)}_{a,\epsilon}(g)\bigr\vert
\ge\delta\epsilon^{-{(2-p)/ d}} \Bigr\}.
\end{eqnarray*}
Applying (\ref{poisson-9}) on the left-hand side,
\begin{eqnarray*}
-I_D(\gamma+\delta)
&\le&\max\Bigl\{\liminf_{\epsilon\to0^+}\epsilon^{2/(d-p)}\log
\P\Bigl\{\inf_{g\in\mathcal{G}_d(D)}G^{(0)}_{a,\epsilon}(g)
\le-\gamma\epsilon^{-{(2-p)/d}} \Bigr\},\\
&&\hspace*{32pt}\limsup_{\epsilon\to0^+}\epsilon^{2/ (d-p)}\log
\P\Bigl\{\sup_{g\in\mathcal{G}_d(D)}\bigl\vert F^{(0)}_{a,\epsilon}(g)\bigr\vert
\ge\delta\epsilon^{-{(2-p)/ d}} \Bigr\} \Bigr\}.
\end{eqnarray*}
Let $a\to\infty$ on the right-hand side. By (\ref{poisson-7}),
\[
\liminf_{a\to\infty}\liminf_{\epsilon\to0^+}\epsilon^{2/ (d-p)}\log
\P\Bigl\{\inf_{g\in\mathcal{G}_d(D)}G^{(0)}_{a,\epsilon}(g)
\le-\gamma\epsilon^{-{(2-p)/d}} \Bigr\}\ge-I_D(\gamma+\delta).
\]
Letting $\delta\to0^+$ on the right-hand side leads to (\ref{poisson-8}).

\subsection{\texorpdfstring{Proof of Theorem \protect\ref{poisson-14}}{Proof of Theorem 3.2}}

Based on (\ref{poisson-23}),
assertion (\ref{poisson-15}) follows from the same argument used in
(\ref{poisson-7}).

By the decomposition
\begin{eqnarray*}
G^{(1)}_{a,\epsilon}(g)
&=&\int_{\R^d} \biggl[\int_DK^{(1)}_{a,\epsilon}(y-x)g^2(y)\,dy \biggr]
\omega(\epsilon\, dx)\\
&&{}-\epsilon\int_{\R^d} \biggl[\int_D
K^{(1)}_{a,\epsilon}(y-x)g^2(y)\,dy \biggr]\,dx,
\end{eqnarray*}
by the uniform [over $g\in\mathcal{G}_d(D)$] bound
\begin{eqnarray*}
\int_{\R^d} \biggl[\int_D
K^{(1)}_{a,\epsilon}(y-x)g^2(y)\,dy \biggr]\,dx
&=& \biggl(\int_Dg^2(y)\,dy \biggr) \biggl(\int_{\R^d}K^{(1)}_{a,\epsilon}(x)\,dx
\biggr)\\
&=&O \biggl( \biggl(\log{1\over\epsilon} \biggr)^{(d-p)/p} \biggr)
\end{eqnarray*}
and by (\ref{poisson-15}), all we need is to establish that
%
\begin{equation}\label{poisson-32}\qquad
\lim_{\epsilon\to0^+}{1\over l(\epsilon)}
\log\P\Bigl\{\sup_{g\in\mathcal{G}_d(D)}\eta_{a,\epsilon}(g)
\ge\gamma\epsilon^{-{(2-p)/ d}} \Bigr\}
= -{2+d-p\over d\rho^*_D}\gamma,
\end{equation}
where
\[
\eta_{a,\epsilon}(g)=\int_{\R^d}
\biggl[\int_DK^{(1)}_{a,\epsilon}(y-x)g^2(y)\,dy \biggr]
\omega(\epsilon\,dx).
\]

Since $K^{(1)}_{a,\epsilon}(y-x)=0$ as
$\vert y-x\vert> 3a(\log\epsilon^{-1})^{1/p}$, $x$ is limited to
a ball with the center 0 and the radius $C(\log\epsilon^{-1})^{1/p}$
when $y\in D$. Consequently,
\begin{eqnarray}\label{poisson-34}
\qquad \sup_{g\in\mathcal{G}_d(D)}
\eta_{a,\epsilon}(g)&\le&\sup_{g\in\mathcal{G}_d(D)}
\int_{\{\vert x\vert\le C(\log\epsilon^{-1})^{1/p}\}}
\biggl[\int_D{g^2(y)\over\vert y-x\vert^p}\,dy \biggr]\omega(\epsilon
\,dx)\nonumber\\
&\le&\rho^*_D\omega\{\vert x\vert\le
C\epsilon^{1/d}
(\log\epsilon^{-1})^{1/p} \}\\
&=&\rho^*_D\widetilde{Z}_\epsilon,\nonumber
\end{eqnarray}
where $\widetilde{Z}_\epsilon\equiv\omega\{\vert x\vert\le
C\epsilon^{1/d}(\log\epsilon^{-1})^{1/p} \}$ is a
Poisson random variable with
\[
\E\widetilde{Z}_\epsilon=\omega_d C^d\epsilon(\log\epsilon^{-1})^{d/p}.
\]

For any $\theta>0$
\[
\E\exp\biggl\{\theta
\biggl(\log{1\over\epsilon} \biggr)\widetilde{Z}_\epsilon\biggr\}
=\exp\{\omega_d C^d\epsilon(\log\epsilon^{-1})^{d/p}
(e^{\theta\log\epsilon^{-1}}-1 ) \}.
\]
Consequently,
\[
\lim_{\epsilon\to0^+}{1\over l(\epsilon)}\log
\E\exp\biggl\{\theta\biggl(\log{1\over\epsilon} \biggr)\widetilde{Z}_\epsilon\biggr\}
=0, \qquad \theta<{2+d-p\over d}.
\]
A standard application of Chebyshev's inequality gives
\[
\limsup_{\epsilon\to0^+}{1\over l(\epsilon)}\log\P\bigl\{\widetilde
{Z}_\epsilon\ge
\gamma\epsilon^{(2-p)/d} \bigr\}\le-{2+d-p\over d}\gamma
\]
for every $\gamma>0$. Thus, the upper bound of
(\ref{poisson-32}) follows from (\ref{poisson-34}).

On the other hand, let $x_0\in\R^d$ be fixed but arbitrary, and
write $\omega_{x_0}(\epsilon \,dx)=\omega(\epsilon(x_0+dx) )$.
Given $\delta>0$ and $\lambda>1$, by variable shifting
%
\begin{eqnarray*}
\hspace*{-1.5pt}&&\sup_{g\in\mathcal{G}_d(D)}\eta_{a,\epsilon}(g)\\
\hspace*{-1.5pt}&&\qquad =\sup_{g\in\mathcal{G}_d(D)}\int_{\R^d}
\biggl[\int_DK^{(1)}_{a,\epsilon}(y-x_0-x)g^2(y)\,dy \biggr]
\omega_{x_0}(\epsilon\, dx) \\
\hspace*{-1.5pt}&&\qquad \ge\sup_{g\in\mathcal{G}_d(D)}\int_{\{\vert x\vert\le\delta\}}
\biggl[\int_DK^{(1)}_{a,\epsilon}(y-x_0-x)g^2(y)\,dy \biggr]
\omega_{x_0}(\epsilon \,dx)\\
\hspace*{-1.5pt}&&\qquad \ge \biggl(\sup_{g\in\mathcal{G}_d(D)}\int_D{1\over(\vert y-x_0\vert
+\delta)^p}
\alpha\bigl(a^{-1}(\log\epsilon^{-1})^{-1/p}(\vert y-x_0\vert+\delta) \bigr)
g^2(y)\,dy \biggr)\\
\hspace*{-1.5pt}&&\qquad \quad
{}\times\omega\{\vert x+x_0\vert\le\epsilon^{1/d}\delta\}
\\
\hspace*{-1.5pt}&&\qquad \stackrel{d}{=} \biggl(\sup_{g\in\mathcal{G}_d(D)}\int_D{g^2(y)\over
(\vert y-x_0\vert+\delta)^p}
\,dy \biggr)\omega\{\vert x\vert\le\epsilon^{1/d}\delta\}
\end{eqnarray*}
as $\epsilon$ is sufficiently small.

Write $Z_\epsilon=\omega\{\vert x\vert\le\epsilon^{1/d}\delta\}$
and $k(\epsilon)= [\gamma\epsilon^{-{(2-p)/d}} ]+1$.
\[
\P\bigl\{Z_\epsilon\ge\gamma\epsilon^{-{(2-p)/d}} \bigr\}
\ge\P\{Z_\epsilon=k(\epsilon) \}=e^{-\omega_d\epsilon\delta^d}
{(\omega_d\epsilon\delta^d)^{k(\epsilon)}\over k(\epsilon)!}.
\]
By Stirling's formula, one can show that for any $\gamma>0$,
\[
\liminf_{\epsilon\to\infty}{1\over l(\epsilon)}\log
\P\bigl\{Z_\epsilon\ge\gamma\epsilon^{-{(2-p)/d}} \bigr\}
\ge-{2+d-p\over d}\gamma,\qquad  \gamma>0.
\]
Replacing $\gamma$ by
\[
\gamma\biggl(\sup_{g\in\mathcal{G}_d(D)}\int_D{g^2(y)
\over(\vert y-x_0\vert+\delta)^p}
\,dy \biggr)^{-1},
\]
we have
\begin{eqnarray*}
&&\liminf_{\epsilon\to\infty}{1\over l(\epsilon)}\log
\P\Bigl\{\sup_{g\in\mathcal{G}_d(D)}\eta_{a,\epsilon}(g)
\ge\gamma\epsilon^{-{(2-p)/d}} \Bigr\}\\
&&\qquad \ge -{2+d-p\over d} \biggl(\sup_{g\in\mathcal{G}_d(D)}
\int_D{g^2(y)\over(\vert y-x_0\vert+\delta)^p}
\,dy \biggr)^{-1}\gamma.
\end{eqnarray*}
Letting $\delta\to0^+$ and taking $x_0\in\R^d$ on the right-hand
side lead to the lower bound of
(\ref{poisson-32}).

\section{Bridging to the eigenvalue problem}\label{FK}

$\!\!\!$Throughout this section, let \mbox{$D\!\subset\!\R^d$} be a bounded open domain,
and let
%
\begin{equation}\label{FK-0}
\mathcal{F}_d(D)= \biggl\{g\in W^{1,2}(D); \int_Dg^2(x)\,dx=1 \biggr\}.
\end{equation}
Given a measurable function $\xi(x)$ on $\R^d$, we introduce the notation
\[
\lambda_\xi(D)=\sup_{g\in\mathcal{F}_d(D)} \biggl\{\int_D\xi
(x)g^2(x)\,dx-{1\over2}
\int_D\vert\nabla g(x)\vert^2\,dx \biggr\}.
\]
Clearly, $\lambda_\xi(D)\le\lambda_\eta(D)$ whenever $\xi(x)\le\eta(x)$
($x\in D$).

Write
\[
\tau_D=\inf\{s\ge0; B_s\notin D\}.
\]

It is well known that by the Feynman--Kac formula,
\[
\E_0 \biggl[\exp\biggl\{\int_0^t\xi(B_s)\,ds \biggr\}; \tau_D\ge t \biggr]
\approx
\exp\{t\lambda_\xi(D) \}\qquad  (t\to\infty)
\]
in some proper sense. For the applications to our setting,
some more explicit bounds are needed. This is our objective in this section.

\begin{lemma}\label{FK-1}
The inequality
%
\begin{equation}\label{FK-2}
\int_D\E_x \biggl[
\exp\biggl\{\int_0^t \xi(B_s)\,ds\biggr \};
\tau_D\ge t \biggr]\,dx\le\vert D\vert\exp\{t\lambda_\xi(D) \}
\end{equation}
holds regardless whether $\lambda_\xi(D)$ is finite or infinite.
\end{lemma}

\begin{pf}
The argument in the case when $\xi(x)\le N$ for some constant $N>0$
is classic (see the treatment given e.g., in \cite{Chen}, Section 4.1):
A standard argument through\vadjust{\goodbreak} a spectral theory [the boundedness of $\xi
(\cdot)$
guarantees the boundedness of the underlined
linear operators in the argument]
gives
that for any $g\in W^{1,2}(D)$
\[
\int_Dg(x)\E_x \biggl[\exp\biggl\{\int_0^t\xi(B_s) \biggr\}g(B_t);
\tau_D\ge t \biggr]\,dx\le\|g\|_{\mathcal{L}^2(D)}^2\exp\{t\lambda_{\xi
}(D) \}.
\]
In particular, let $g_n\in W^{1,2}(D)$ be a monotonic
sequence such that $0\le g_n(x)\le1$ and $g_n(x)\uparrow1$ ($n\to
\infty$) for
every $x\in D$. Then
\begin{eqnarray}
\int_Dg_n(x)\E_x \biggl[\exp\biggl\{\int_0^t\xi(B_s) \biggr\}g_n(B_t);
\tau_D\ge t \biggr]\,dx\le\vert D\vert\exp\{t\lambda_{\xi}(D) \},\nonumber
\\
\eqntext{n=1,2,\ldots.}
\end{eqnarray}
Letting $n\to\infty$ on the left-hand side, the desired bound follows
from monotonic convergence.

To remove the boundedness assumption, we write $\xi_N(x)=\min\{\xi(x),
N\}$.
By what has been proved,
\[
\int_D\E_x \biggl[\exp\biggl\{\int_0^t\xi_N(B_s)\,ds \biggr\};
\tau_D\ge t \biggr]\,dx\le\vert D\vert\exp\{t
\lambda_{\xi_N}(D) \}\le\vert D\vert\exp\{t
\lambda_{\xi}(D) \}.
\]
The conclusion follows from monotonic convergence again
as we let $N\to\infty$
on the left-hand side.
\end{pf}

\begin{lemma}\label{FK-1'}
For any $\alpha,\beta>1$ satisfying
$\alpha^{-1}+\beta^{-1}=1$ and $\lambda_{(\beta/\alpha)\xi}(D)<\infty$
[in this case $\lambda_{\alpha^{-1}\xi}(D)<\infty$] and $0<\delta<t$
\begin{eqnarray}\label{FK-2'}
&&\int_D\E_x \biggl[
\exp\biggl\{\int_0^t \xi(B_s)\,ds \biggr\};
\tau_D\ge t \biggr]\,dx\nonumber\\
&&\qquad\ge(2\pi)^{\alpha d/2}\delta^{d/2}t^{\alpha d/ (2\beta)}
\vert D\vert^{-2\alpha/\beta}\\
&&\qquad \quad
{}\times\exp\bigl\{-\delta(\alpha/\beta)
\lambda_{(\beta/\alpha)\xi}(D) \bigr\}
\exp\{\alpha(t+\delta)\lambda_{\alpha^{-1}\xi}(D) \}.\nonumber
\end{eqnarray}
\end{lemma}

\begin{pf}
We only need to show that
%
\begin{eqnarray}\label{FK-3'}
&&\int_D\E_x \biggl[
\exp\biggl\{\int_0^t \xi(B_s)\,ds \biggr\};
\tau_D\ge t \biggr]\,dx\nonumber\\
&&\qquad\ge(2\pi)^{\alpha d/2}\delta^{d/2}t^{\alpha d/ (2\beta)}
\vert D\vert^{-\alpha/\beta}
\exp\{\alpha(t+\delta)\lambda_{\alpha^{-1}\xi}(D) \}\\
&&\qquad\quad{}\times\biggl\{\int_D\E_x \biggl[
\exp\biggl\{{\beta\over\alpha}\int_0^\delta\xi(B_s)\,ds \biggr\};
\tau_D\ge\delta\biggr]\,dx \biggr\}^{-\alpha/\beta}\nonumber
\end{eqnarray}
as, by Lemma \ref{FK-1},
\[
\int_D\E_x \biggl[
\exp\biggl\{{\beta\over\alpha}\int_0^\delta\xi(B_s)\,ds \biggr\};
\tau_D\ge\delta\biggr]\,dx\le\vert D\vert\exp\bigl\{-\delta(\alpha/\beta)
\lambda_{(\beta/\alpha)\xi}(D) \bigr\}.
\]

We first consider the case when
$\xi(x)$ is H\"{o}lder continuous. By the Feynman--Kac representation,
\[
u(t,x)=\E_x \biggl[
\exp\biggl\{\int_0^t \xi(B_s)\,ds \biggr\};
\tau_D\ge t \biggr]
\]
solves the initial-boundary value problem
\[
\cases{
\partial_tu(t,x)={1\over2}\Delta u(t,x)+\xi(x)u(t,x),&\quad $ (t,x)\in(0, t)
\times D$,\cr
u(0, x)=1,&\quad $ x\in D $,\cr
u(t, x)=0,&\quad $ (t,x)\in(0,\infty)\times\partial D$.
}
\]

Let $\lambda_1>\lambda_2\ge\lambda_3\ge\cdots$ be
the eigenvalues of the operator $(1/2)\Delta+\xi$ in $\mathcal{L}^2(D)$
with zero boundary condition and initial value 1 in $D$, and let
$e_k\in\mathcal{L}^2(D)$
be an orthonormal basis corresponding to $\{\lambda_k\}$. By (2.31)
in~\cite{G-K-M},
\[
\E_x\biggl[\exp\biggl\{\int_0^t \xi(B_s)\,ds \biggr\}\delta_x(B_t);
\tau_D\ge t \biggr]=\sum_{k=1}^\infty e^{t\lambda_k} e^2_k(x)
\ge e^{t\lambda_1} e^2_1(x).
\]
Noticing the fact that $\lambda_1=\lambda_\xi(D)$ and integrating both sides
we have
\[
\int_D\E_x\biggl[\exp\biggl\{\int_0^t \xi(B_s)\,ds \biggr\}\delta_x(B_t);
\tau_D\ge t \biggr]\,dx\ge\exp\{t\lambda_\xi(D) \}.
\]

Replace $\xi$ by $\alpha^{-1}\xi$ and $t$ by $t+\delta$.
By H\"{o}lder's inequality,
\begin{eqnarray*}
&&\exp\{(t+\delta)\lambda_{\alpha^{-1}\xi}(D) \}\\
&&\qquad\le
\int_D\E_x \biggl[
\exp\biggl\{\alpha^{-1}\int_0^{t+\delta} \xi(B_s)\,ds \biggr\}
\delta_x(B_{t+\delta});
\tau_D\ge t+\delta\biggr]\,dx\\
&&\qquad\le\biggl\{\int_D\E_x \biggl[
\exp\biggl\{(\beta/\alpha)\int_t^{t+\delta} \xi(B_s)\,ds \biggr\};
\tau_D\ge t+\delta\biggr]\,dx \biggr\}^{1/\beta}\\
&&\quad\qquad{}\times\biggl\{\int_D\E_x \biggl[\exp\biggl\{\int_0^{t} \xi(B_s)\,ds \biggr\}
\delta_x(B_{t+\delta});
\tau_D\ge t+\delta\biggr]\,dx \biggr\}^{1/\alpha}.
\end{eqnarray*}
Notice that
\begin{eqnarray*}
&&\E_x \biggl[\exp\biggl\{\int_0^{t} \xi(B_s)\,ds \biggr\}\delta_x(B_{t+\delta});
\tau_D\ge t+\delta\biggr]\\
&&\qquad\le\E_x \biggl[\exp\biggl\{\int_0^{t} \xi(B_s)\,ds \biggr\}\delta_x(B_{t+\delta});
\tau_D\ge t \biggr]\\
&&\qquad=\E_x \biggl[\exp\biggl\{\int_0^{t} \xi(B_s)\,ds \biggr\}p_\delta(B_t-x);
\tau_D\ge t \biggr],
\end{eqnarray*}
where
\[
p_\delta(y)={1\over(2\pi\delta)^{d/2}}
\exp\biggl\{-{\vert y\vert^2\over2\delta} \biggr\}\le{1\over(2\pi\delta)^{d/2}}.\vadjust{\goodbreak}
\]

In addition,
\begin{eqnarray*}
&&\int_D\E_x \biggl[
\exp\biggl\{(\beta/\alpha)\int_t^{t+\delta} \xi(B_s)\,ds \biggr\};
\tau_D\ge t+\delta\biggr]\,dx\\
&&\qquad \le\int_D\E_x \biggl[
\exp\biggl\{(\beta/\alpha)\int_t^{t+\delta} \xi(B_s)\,ds \biggr\};
B_t\in D, \tau_D'\ge t+\delta\biggr]\,dx\cr
&&\qquad =\int_D \biggl[\int_Dp_t(y-x)\E_y \biggl(\exp\biggl\{(\beta/\alpha)\int_0^{\delta}
\xi(B_s)\,ds \biggr\}; \tau_D\ge\delta\biggr)\,dy \biggr]\,dx\\
&&\qquad \le{1\over(2\pi t)^{d/2}}\vert D\vert
\int_D\E_y \biggl[\exp\biggl\{(\beta/\alpha)\int_0^{\delta}
\xi(B_s)\,ds \biggr\}; \tau_D\ge\delta\biggr]\,dy,
\end{eqnarray*}
where
\[
\tau_D'=\inf\{s\ge t; B_s\notin D\}.
\]

Summarizing our argument, we have established the bound (\ref{FK-3'}).

We now move to the
case when $\xi(x)\ge-N$ for some $N>0$. For any H\"{o}lder-continuous
$\eta(x)$ on $D$ with $\eta(x)\le\xi(x)$ a.e. on $D$,
\begin{eqnarray*}
&&\int_D\E_x \biggl[\exp\biggl\{\int_0^t\xi(B_s)\,ds \biggr\};
\tau_D\ge t \biggr]\,dx\\
&&\qquad \ge(2\pi)^{\alpha d/2}\delta^{d/2}t^{\alpha d/(2\beta)}
\vert D\vert^{-\alpha/\beta}
\exp\{\alpha(t+\delta)\lambda_{\alpha^{-1}\eta}(D) \}\\
&&\qquad\quad  {}\times\biggl\{\int_D\E_x \biggl[
\exp\biggl\{{\beta\over\alpha}\int_0^\delta\eta(B_s)\,ds \biggr\};
\tau_D\ge\delta\biggr]\,dx \biggr\}^{-\alpha/\beta}\\
&&\qquad \ge(2\pi)^{\alpha d/2}\delta^{d/2}t^{\alpha d/(2\beta)}
\vert D\vert^{-\alpha/\beta}
\exp\{\alpha(t+\delta)\lambda_{\alpha^{-1}\eta}(D) \}\\
&&\qquad \quad {}\times\biggl\{\int_D\E_x \biggl[
\exp\biggl\{{\beta\over\alpha}\int_0^\delta\xi(B_s)\,ds \biggr\};
\tau_D\ge\delta\biggr]\,dx \biggr\}^{-\alpha/\beta}.
\end{eqnarray*}
Let
\[
\mathcal{H}_\xi= \{\eta(\cdot); \mbox{$\eta(x)$ is H\"{o}lder continuous
on $D$ and $\eta(x)\le\xi(x)$ a.e. on $D$} \}.
\]
Since $\xi(\cdot)\ge-N$, $\mathcal{H}_\xi\not=\phi$. Further,
by standard approximation theory, $\mathcal{H}_\xi$ is rich enough to
approximate $\xi$. More precisely, the desired bound follows from
\begin{eqnarray*}
\sup_{\eta\in\mathcal{H}_\xi} \lambda_{\alpha^{-1}\eta}(D)
&=&\sup_{g\in\mathcal{F}_d(D)} \biggl\{\alpha^{-1}\sup_{\eta\in\mathcal{H}_\xi}
\int_D\eta(x)g^2(x)\,dx
-{1\over2}
\int_D\vert\nabla g(x)\vert^2\,dx \biggr\}\\
&=&\lambda_{\alpha^{-1}\xi}(D).
\end{eqnarray*}

To remove the boundedness assumption, we write $\xi_N(x)=\xi(x)\vee(-N)$.
We have
\begin{eqnarray*}
&&\int_D\E_x \biggl[\exp\biggl\{\int_0^t\xi_N(B_s)\,ds \biggr\};
\tau_D\ge t \biggr]\,dx\\
&&\qquad \ge(2\pi)^{\alpha d/2}\delta^{d/2}t^{\alpha d/(2\beta)}
\vert D\vert^{-\alpha/\beta}
\exp\{\alpha(t+\delta)\lambda_{\alpha^{-1}\xi_N}(D) \}\\
&&\qquad \quad {}\times\biggl\{\int_D\E_x \biggl[
\exp\biggl\{{\beta\over\alpha}\int_0^\delta\xi_N(B_s)\,ds \biggr\};
\tau_D\ge\delta\biggr]\,dx \biggr\}^{-\alpha/\beta}.
\end{eqnarray*}
Noticing $\lambda_{\alpha^{-1}\xi_N}(D)\ge\lambda_{\alpha^{-1}\xi}(D)$ and
letting $N\to\infty$,
the monotonic convergence theorem leads to (\ref{FK-3'}).
\end{pf}

\begin{lemma}\label{FK-3} Let $0<\delta<t$, and assume $0\in D$.
%
\begin{eqnarray}\label{FK-4}
&&\E_0 \biggl[\exp\biggl\{\int_0^t \xi(B_s)\,ds \biggr\};
\tau_D\ge t \biggr]\hspace*{-35pt}\nonumber\\
&&\qquad\le\biggl(\E_0\exp\biggl\{\beta\int_0^\delta
\xi(B_s)\,ds \biggr\} \biggr)^{1/\beta}\hspace*{-35pt}\\
&&\qquad\quad{}\times\biggl\{{1\over(2\pi\delta)^{d/2}}
\int_D\E_x \biggl[
\exp\biggl\{\alpha\int_0^{t-\delta} \xi(B_s)\,ds \biggr\};
\tau_D\ge t-\delta\biggr]\,dx \biggr\}^{1/\alpha}.\nonumber\hspace*{-35pt}
\end{eqnarray}

On the other hand,
%
\begin{eqnarray}\label{FK-5}
&&\E_0\exp\biggl\{\int_0^t \xi(B_s)\,ds \biggr\}\nonumber\hspace*{-35pt}\\
&&\qquad\ge\biggl(\E_0\exp\biggl\{-{\beta\over\alpha}\int_0^\delta
\xi(B_s)\,ds \biggr\} \biggr)^{-\alpha/\beta}\hspace*{-35pt}\\
&&\qquad\quad{}\times\biggl\{
\int_Dp_\delta(x)\E_x \biggl[
\exp\biggl\{\alpha^{-1}\int_0^{t-\delta} \xi(B_s)\,ds \biggr\};
\tau_D\ge t-\delta\biggr]\,dx \biggr\}^\alpha,\nonumber\hspace*{-35pt}
\end{eqnarray}
where $p_\delta(x)$ is the density of $B_\delta$.
\end{lemma}

\begin{pf}
By H\"{o}lder's inequality,
\begin{eqnarray*}
&&\E_0 \biggl[\exp\biggl\{\int_0^t \xi(B_s)\,ds \biggr\};
\tau_D\ge t \biggr]\\
&&\qquad \le\biggl(\E_0\exp\biggl\{\beta\int_0^\delta
\xi(B_s)\,ds \biggr\} \biggr)^{1/\beta}
\biggl\{\E_0 \biggl[\exp\biggl\{\alpha\int_\delta^t \xi(B_s)\,ds \biggr\};
\tau_D\ge t \biggr] \biggr\}^{1/\alpha}.
\end{eqnarray*}

Write $\tau_D'=\inf\{s\ge\delta; B_s\notin D\}$.
Claim (\ref{FK-4}) follows from the following procedure via
Markov property:
\begin{eqnarray*}
&&\E_0 \biggl[\exp\biggl\{\alpha\int_\delta^t \xi(B_s)\,ds \biggr\};
\tau_D\ge t \biggr]\\
&&\qquad \le\E_0 \biggl[\exp\biggl\{\alpha\int_\delta^t \xi(B_s)\,ds \biggr\};
B_\delta\in D, \tau_D'\ge t \biggr]\\
&&\qquad =\int_Dp_\delta(x)\E_x \biggl[\exp\biggl\{\alpha\int_0^{t-\delta} \xi(B_s)\,ds \biggr\};
\tau_D\ge t-\delta\biggr]\,dx\\
&&\qquad \le{1\over(2\pi\delta)^{d/2}}\int_D
\E_x \biggl[\exp\biggl\{\alpha\int_0^{t-\delta} \xi(B_s)\,ds \biggr\};
\tau_D\ge t-\delta\biggr]\,dx.
\end{eqnarray*}

On the other hand,
\begin{eqnarray*}
&&\E_0 \biggl[\exp\biggl\{\alpha^{-1}\int_\delta^t \xi(B_s)\,ds \biggr\};
B_\delta\in D, \tau_D'\ge t \biggr]\\
&&\qquad \le\E_0 \biggl[\exp\biggl\{-\alpha^{-1}
\int_0^\delta\xi(B_s)\,ds \biggr\}
\exp\biggl\{\alpha^{-1}\int_0^t \xi(B_s)\,ds \biggr\} \biggr]\\
&&\qquad \le\biggl(\E_0\exp\biggl\{-{\beta\over\alpha}\int_0^\delta
\xi(B_s)\,ds \biggr\} \biggr)^{1/\beta} \biggl\{\E_0
\exp\biggl\{\int_0^{t} \xi(B_s)\,ds \biggr\} \biggr\}^{1/\alpha}.
\end{eqnarray*}
Thus, (\ref{FK-5}) follows from Markov property which claims that
\begin{eqnarray*}
&&\E_0 \biggl[\exp\biggl\{\alpha^{-1}\int_\delta^t \xi(B_s)\,ds \biggr\};
B_\delta\in D,
\tau_D'\ge t \biggr]\\
&&\qquad=\int_Dp_\delta(x)\E_x \biggl[\exp\biggl\{\alpha^{-1}\int
_0^{t-\delta}
\xi(B_s)\,ds \biggr\};
\tau_D\ge t-\delta\biggr]\,dx.
\end{eqnarray*}
\upqed
\end{pf}

\section{Upper bounds}\label{u}

In this section we establish the upper bounds
for Theorems~\ref{th-2} and \ref{th-3}.
More precisely, we prove that
%
\begin{eqnarray}\label{u-1}
&&\limsup_{t\to\infty}t^{-1}(\log t)^{-{(d-p)/ d}}\nonumber\\[-8pt]\\[-8pt]
&&\hphantom{\limsup_{t\to\infty}} {}\times\log
\E_0\exp\biggl\{-\theta\int_0^t
\ol{V}(B_s)\,ds \biggr\}\le\Lambda_0(\theta) \qquad\mbox{a.s.-}\P\nonumber
\end{eqnarray}
when $d/2<p<d$, and
%
\begin{eqnarray}\label{u-2}
&&\limsup_{t\to\infty}{1\over t} \biggl({\log\log t\over\log t}
\biggr)^{2/(2-p)}\nonumber\\[-8pt]\\[-8pt]
&&\hphantom{\limsup_{t\to\infty}} {}\times
\log\E_0\exp\biggl\{\theta\int_0^t\ol{V}(B_s)\,ds \biggr\}\le\Lambda
_1(\theta)\qquad \mbox{a.s.-}\P\nonumber
\end{eqnarray}
when $d/2<p<\min\{2, d\}$,
where
%
\begin{eqnarray}\label{u-1'}
\Lambda_0(\theta)&=&{\theta d^2\over d-p}
\biggl({\omega_d\over d}\Gamma\biggl({2p-d\over p} \biggr) \biggr)^{p/d},
\\
\label{u-2'}
\Lambda_1(\theta)&=&{1\over2}p^{p/(2-p)}
(2-p)^{(4-p)/(2-p)} \biggl({d\theta\,\sigma(d,p)\over2+d-p} \biggr)^{2/(2-p)}.
\end{eqnarray}

The following notation\vspace*{1pt} will be used in this and the next sections.
For any \mbox{$R>0$}, $Q_R=(-R, R)^d$.
%
\begin{equation}\label{u-3}
h_t= \cases{
(\log t)^{(d-p)/(2d)}, &\quad for the proof of (\ref{u-1}),\cr
\displaystyle \biggl({\log t\over\log\log t} \biggr)^{1/(2-p)},&\quad
for the proof of (\ref{u-2}).
}
\end{equation}

Write $R_k=R_k(t)= (Mth_t)^k$ ($k=1,2,\ldots$) where the constant $M>0$
is fixed but sufficiently large. Write $\xi(x)=-\ol{V}(x)$ in the proof of
(\ref{u-1}) and $\xi(x)=\ol{V}(x)$ in the proof of
(\ref{u-2}).\vspace*{-1pt}

Finally we recall that for any open domain $D\subset\R^d$ containing 0,
\[
\tau_D=\inf\{s\ge0; B_s\notin D\}.
\]

Consider
the decomposition
%
\begin{eqnarray}\label{u-4}
&&\E_0\exp\biggl\{\theta\int_0^t\xi(B_s)\,ds \biggr\}\nonumber\\
&&\qquad=\E_0 \biggl[\exp\biggl\{\theta\int_0^t\xi(B_s)\,ds \biggr\};
\tau_{Q_{R_1}}\ge t \biggr]\nonumber\\
&&\quad\qquad{}+\sum_{k=1}^\infty\E_0 \biggl[\exp\biggl\{\theta\int_0^t\xi(B_s)\,
ds \biggr\};
\tau_{Q_{R_k}}< t\le\tau_{Q_{R_k}} \biggr]\nonumber\\
&&\qquad\le\E_0 \biggl[\exp\biggl\{\theta\int_0^t\xi(B_s)\,ds \biggr\};
\tau_{Q_{R_1}}\ge t \biggr]\\
&&\qquad\quad{}+\sum_{k=1}^\infty(\P\{\tau_{Q_{R_k}}< t\} )^{1/2}\nonumber
\\
&&\qquad \quad \hphantom{{}+\sum_{k=1}^\infty}
{}\times\biggl\{\E_0 \biggl[\exp\biggl\{2\theta\int_0^t\xi(B_s)\,ds \biggr\};
\tau_{Q_{R_{k+1}}}\ge t \biggr] \biggr\}^{1/2}.\nonumber
\end{eqnarray}

The well-known result on the Gaussian tail gives that
\[
(\P\{\tau_{Q_{R_k}}< t\} )^{1/2}\le\exp\{-c R_k^2/t \}
=\exp\{-c M^2t^{2k-1}h_t^{2k} \}.
\]

Let $\alpha,\beta>1$ satisfy $\alpha^{-1}+\beta^{-1}=1$ with $\alpha$
close to 1. By (\ref{FK-4}) (with $\delta=1$) and Lemma \ref{FK-1},
\begin{eqnarray*}
&&\E_0 \biggl[\exp\biggl\{\theta\int_0^t\xi(B_s)\,ds \biggr\};
\tau_{Q_{R_1}}\ge t \biggr]\\
&&\qquad\le
{1\over(2\pi)^{d/\alpha}}
\biggl(\E_0\exp\biggl\{\theta\beta\int_0^{1}\xi_{R_1}(B_s)\,ds \biggr\}
\biggr)^{1/\beta}\\
&&\quad\qquad{}\times\biggl\{\int_{Q_{R_1}}\,dx\,
\E_x \biggl[\exp\biggl\{\theta\alpha\int_0^1\xi(B_s)\,ds \biggr\};
\tau_{Q_{R_1}}\ge t-1 \biggr] \biggr\}^{1/\alpha}\\
&&\qquad\le\biggl({R_1\over\pi} \biggr)^{d/\alpha}
\biggl(\E_0\exp\biggl\{\theta\beta\int_0^1\xi(B_s)\,ds \biggr\} \biggr)^{1/\beta}
\exp\{(t-1)\lambda_{\theta\alpha\xi}(Q_{R_1}) \}.
\end{eqnarray*}
Similarly,
\begin{eqnarray*}
&&\E_0 \biggl[\exp\biggl\{2\theta\int_0^t\xi(B_s)\,ds \biggr\};
\tau_{\widetilde{Q}_{R_{k+1}}}\ge t \biggr]\\
&&\qquad\le\biggl({R_{k+1}\over\pi} \biggr)^{d/\alpha}\!
\biggl(\E_0\exp\biggl\{2\theta\beta\int_0^1\xi(B_s)\,ds \biggr\} \biggr)^{1/\beta}\!
\exp\{(t-1)\lambda_{2\theta\alpha\xi}(Q_{R_{k+1}}) \}.
\end{eqnarray*}

Summarizing our estimates since (\ref{u-4}),
%
\begin{eqnarray}\label{u-5}
&&\E_0\exp\biggl\{\theta\int_0^t\xi(B_s)\,ds \biggr\}\nonumber\hspace*{-30pt}\\
&&\qquad\le\biggl({R_1\over\pi} \biggr)^{d/\alpha}
\biggl(\E_0\exp\biggl\{\theta\beta\int_0^1\xi(B_s)\,ds \biggr\} \biggr)^{1/\beta}
\exp\{t\lambda_{\alpha\theta\xi}(Q_{R_1}) \}\nonumber\hspace*{-30pt}\\[-8pt]\\[-8pt]
&&\qquad\quad{}+ \biggl(\E_0\exp\biggl\{2\theta\beta\int_0^1\xi(B_s)\,ds \biggr\}
\biggr)^{1/2\beta}\nonumber\hspace*{-30pt}\\
&&\qquad\qquad{}\times\sum_{k=1}^\infty
\biggl({R_{k+1}\over\pi} \biggr)^{d/2\alpha}\exp\{-c M^2t^{2k-1}h_t^{2k} \}\exp\biggl\{
{t\over2}
\lambda_{2\alpha\theta\xi}(Q_{R_{k+1}}) \biggr\}.\nonumber\hspace*{-30pt}
\end{eqnarray}

To prove (\ref{u-1}) and (\ref{u-2}), therefore, all we need is to show that
for any $\theta>0$,
%
\begin{eqnarray}\label{u-6}
\lim_{t\to\infty}h_t^{-2}\lambda_{\theta\xi}(Q_t)
\le\Lambda(\theta)\equiv\cases{
\Lambda_0(\theta),&\quad
for the proof of (\ref{u-1}),\cr
\Lambda_1(\theta),&\quad
for the proof of (\ref{u-2}).
}\nonumber\\[-8pt]\\[-8pt]
\eqntext{\mbox{a.s.-}\P.}
\end{eqnarray}

Indeed, we apply (\ref{u-6}) to the first term on the right-hand side of
(\ref{u-5}) (with $t$ being replaced by $R_1=Mth_t$ and $\theta$ being
replaced by $\alpha\theta$). Notice that $\alpha$ can be arbitrarily close
to 1. This term alone
does not exceed the limit set in (\ref{u-1}) and (\ref{u-2}) if we let
$\alpha\to1^+$ after the limit for $t$.
To control the infinite series on the right-hand side of (\ref{u-5}),
we apply (\ref{u-6}) to each term with $t$ being replaced by
$R_{k+1}=(Mth_t)^{k+1}$ and with $\theta$ being replaced by~$2\alpha
\theta$.
In this way, the series is dominated by
\[
\sum_{k=1}^\infty\biggl({R_{k+1}\over\pi} \biggr)^{d/\alpha}
\exp\{-c't^{2k-2}h_t^{2k} \}
=O(1) \qquad \mbox{a.s.-}\P\ (t\to\infty),
\]
where $c'>0$ is a constant. Here we point out that to control the first term
of the series in (\ref{u-5}), $M>0$ is required to be sufficiently
large.\vadjust{\goodbreak}

Let $\delta>0$ be a small number, and write
\[
\tilde{h}_t=h_t\sqrt{u\over1+\delta}.\vspace*{-2pt}
\]
Define
\[
\xi_t(x)=
\pm\theta\tilde{h}_t^{p-2}\int_{\R^d}{1\over\vert y-x\vert^p}
[\omega(\tilde{h}_t^{-d}dx )-\tilde{h}_t^{-d}\,dx ],\vspace*{-2pt}
\]
where ``$-$'' corresponds to the proof of (\ref{u-1}) and ``$+$''
corresponds to the proof
of (\ref{u-2}).

Under the substitution
\[
g(x)\mapsto\tilde{h}_t^{d/2}g (x\tilde{h}_t ),\vspace*{-2pt}
\]
we have that
\[
\lambda_{\theta\xi}(Q_{t})=\tilde{h}_t^2\sup_{g\in\mathcal
{F}_d(Q_{t\tilde{h}_t})}
\biggl\{\int_{Q_{t\tilde{h}_t}}\xi_t(x)g^2(x)\,dx-{1\over2}
\int_{Q_{t\tilde{h}_t}}\vert\nabla g(x)\vert^2\,dx \biggr\}.\vspace*{-2pt}
\]

Let $r\ge2$ be large but fixed. By Proposition 1 in \cite{G-K}, or by
Lemma 4.6 in \cite{G-K-M}, there is
a nonnegative and continuous function $\Phi(x)$ on $\R^d$
whose support is contained in the 1-neighborhood of the grid $2r\Z^d$,
such that
\[
\lambda_{\xi_t-\Phi^y}(Q_{t\tilde{h}_t})
\le\max_{z\in2r\Z^d\cap Q_{2t\tilde{h}_t+2r}}\lambda_{\xi_t}
(z+Q_{r+1}),\qquad
y\in Q_r,\vspace*{-2pt}
\]
where $\Phi^y(x)=\Phi(x+y)$.
In addition, $\Phi(x)$ is periodic with period $2r$
\[
\Phi(x+2rz)=\Phi(x);\qquad  x\in\R^d, z\in\Z^d,\vspace*{-2pt}
\]
and there is a constant $K>0$ independent of $r$ and $t$ such that
\[
\int_{Q_{r}}\Phi(x)\,dx\le{K\over r}.\vspace*{-2pt}
\]

By periodicity
\begin{eqnarray*}
&&\sup_{g\in\mathcal{F}_d(Q_{t\tilde{h}_t})}
\biggl\{\int_{Q_{t\tilde{h}_t}}\xi_t(x)g^2(x)\,dx-{1\over2}
\int_{Q_{t\tilde{h}_t}}\vert\nabla g(x)\vert^2\,dx \biggr\}\\[-3pt]
&&\qquad \le{K\over r(2r)^d}+\sup_{g\in\mathcal{F}_d(Q_{t\tilde{h}_t})}
\biggl\{\int_{Q_{t\tilde{h}_t}} \biggl(\xi_t(x)-{1\over(2r)^d}
\int_{Q_r}
\Phi^y(x)dy \biggr)g^2(x)\,dx\\[-3pt]
&&\hspace*{224pt}
{}-{1\over2}
\int_{Q_{t\tilde{h}_t}}\vert\nabla g(x)\vert^2\,dx \biggr\}\\[-3pt]
&&\qquad \le{K\over r(2r)^d}+{1\over(2r)^d}\int_{Q_r}\sup_{g\in\mathcal
{F}_d(Q_{t\tilde{h}_t})}
\biggl\{\int_{Q_{t\tilde{h}_t}} \bigl(\xi_t(x)-\Phi^y(x) \bigr)
g^2(x)\,dx\\[-3pt]
&&\hspace*{200pt}
{}-{1\over2}
\int_{Q_{Q_{t\tilde{h}_t}}}\vert\nabla g(x)\vert^2\,dx \biggr\}\,dy\\[-3pt]
&&\qquad ={K\over r(2r)^d}+{1\over(2r)^d}\int_{Q_{r}}\lambda_{\xi_t-\Phi^y}
(Q_{t\tilde{h}_t})\,dy\\[-3pt]
&&\qquad \le{K\over2^dr^{d+1}}+\max_{z\in2r\Z^d\cap Q_{2t\tilde{h}_t+2r}}
\lambda_{\xi_t}(z+Q_{r+1}).
\end{eqnarray*}

Summarizing our estimates
\[
\lambda_{\theta\xi}(Q_t)\le{uh_t^2\over1+\delta} \biggl\{{K\over2^dr^{d+1}}+
\max_{z\in2r\Z^d\cap Q_{2t\tilde{h}_t+2r}}
\lambda_{\xi_t}(z+Q_{r+1}) \biggr\}.
\]
Take $r>0$ sufficiently large so that the first term on the right-hand side
is less than ${\delta u\over1+\delta}h_t^2$. We have that
%
\begin{equation}\label{u-7}
\P\{\lambda_{\theta\xi}(Q_t)\ge uh_t^2 \}\le
\P\Bigl\{\max_{z\in2r\Z^d\cap Q_{2t\tilde{h}_t+2r}}
\lambda_{\xi_t}(z+Q_{r+1})>1 \Bigr\}.
\end{equation}

By shifting invariance of the Poisson field, the random variables
\[
\lambda_{\xi_t}(z+Q_{r+1});\qquad  z\in2r\Z^d\cap Q_{2t\tilde{h}_t+2r},
\]
are identically distributed. Consequently, there is $C>0$
%
\begin{eqnarray}\label{u-8}
&&\P\Bigl\{\max_{z\in2r\Z^d\cap Q_{2t\tilde{h}_t+2r}}
\lambda_{\xi_t}(z+Q_{r+1})>1 \Bigr\}\nonumber \\
&&\qquad \le C(th_t)^d\P\{\lambda_{\xi_t}(Q_{r+1})> 1\}\\
&&\qquad=C(th_t)^d\P\Bigl\{\sup_{g\in\mathcal{G}_d(Q_{r+1})}\int
_{Q_{r+1}}\xi_t(x)g^2(x)\,dx
>1 \Bigr\}\nonumber,
\end{eqnarray}
where $\mathcal{G}_d(Q_{r+1})$ is defined in (\ref{poisson-1}) and the last
step follows from Lemma~\ref{A-1}.

We now reach the point of applying Theorems \ref{poisson-6} and
\ref{poisson-14}. In connection with (\ref{u-1}),
\begin{eqnarray*}
&&\sup_{g\in\mathcal{G}_d(Q_{r+1})}\int_{Q_{r+1}}\xi_t(x)g^2(x)\,dx\\
&&\qquad =-\theta\tilde{h}_t^{p-2}\inf_{g\in\mathcal{G}_d(Q_{r+1})}
\int_{\R^d} \biggl[\int_{Q_{r+1}}{g^2(y)\over\vert y-x\vert^p}\,dy
\biggr] [\omega(\tilde{h}_t^{-d}\,dx)-\tilde{h}_t^{-d}\,dx ].
\end{eqnarray*}
Taking $\epsilon=\tilde{h}_t^{-d}$ and $\gamma=\theta^{-1}$
in (\ref{poisson-9}) leads to
\begin{eqnarray}\label{u-10}
&&\lim_{t\to\infty}{1\over\log t}\log
\P\biggl\{\sup_{g\in\mathcal{G}_d(Q_{r+1})}\int_{Q_{r+1}}\xi_t(x)g^2(x)\,dx
>1 \biggr\}\nonumber\\
&&\qquad=- \biggl({u\over1+\delta} \biggr)^{d/(d-p)}I_{Q_{r+1}}
(\theta^{-1})\\
&&\qquad \le- \biggl({u(d-p)\over\theta d(1+\delta)} \biggr)^{d/(d-p)}
\biggl(\omega_d\Gamma\biggl({2p-d\over p} \biggr) \biggr)^{-{p/(d-p)}}\nonumber,
\end{eqnarray}
where the rate function $I_{Q_{r+1}}(\cdot)$
is defined in (\ref{poisson-10}),
and the last step follows from the obvious fact that
\[
\sup_{g\in\mathcal{G}_d(Q_{r+1})}
\|g\|_{\mathcal{L}^2(Q_{r+1})}\le1.
\]

Take $u=(1+2\delta)\Lambda(\theta)$.
By (\ref{u-7}), (\ref{u-8}) and (\ref{u-10}), there is a $\nu>0$ such that
%
\begin{equation}\label{u-11}
 \P\{\lambda_{\theta\xi}(Q_t)\ge(1+2\delta)\Lambda(\theta)h_t^2 \}
\le C(th_t)^d\exp\{(d+\nu)\log t \}
=C{h_t^d\over t^{\nu}}\hspace*{-35pt}
\end{equation}
for sufficiently large $t$.

We now establish (\ref{u-11}) for the proof of (\ref{u-2}).
In this case
\begin{eqnarray*}
&&\sup_{g\in\mathcal{G}_d(Q_{r+1})}
\int_{Q_{r+1}}\xi_t(x)g^2(x)\,dx\\
&&\qquad =\theta\tilde{h}_t^{p-2}\sup_{g\in\mathcal{G}_d(Q_{r+1})}
\int_{\R^d} \biggl[\int_{Q_{r+1}}{g^2(y)\over\vert y-x\vert^p}\,dy
\biggr] [\omega(\tilde{h}_t^{-d}dx)-\tilde{h}_t^{-d}\,dx ].
\end{eqnarray*}
Taking $\epsilon=\tilde{h}_t^{-d}$ and $\gamma=\theta^{-1}$
in (\ref{poisson-17}),
\begin{eqnarray*}
&&\lim_{t\to\infty}{1\over\log t}\log
\P\biggl\{\sup_{g\in\mathcal{G}_d(Q_{r+1})}\int_{Q_{r+1}}\xi_t(x)g^2(x)\,dx
>1 \biggr\}\\
&&\qquad = - \biggl({u\over1+\delta} \biggr)^{(2-p)/2}
{2+d-p\over\theta(2-p)\rho^*_{Q_{r+1}}},
\end{eqnarray*}
where $\rho^*_{Q_{r+1}}$ is defined as the second variation in
(\ref{poisson-13}) with $D=Q_{r+1}$.

Write
%
\begin{eqnarray}\label{u-11'}
\mathcal{G}_d&=&\mathcal{G}_d(\R^d)= \biggl\{g\in W^{1,2}(\R^d);
\|g\|_2+{1\over2}\|\nabla g\|_2^2=1 \biggr\},
\\
\label{u-11''}
\rho(d,p)&=&\sup_{g\in\mathcal{G}_d}\int_{\R^d}{g^2(x)\over\vert y\vert^p}
\,dy \quad\mbox{and}\nonumber\\[-8pt]\\[-8pt]
\rho^*(d,p)&=&\sup_{g\in\mathcal{G}_d}
\sup_{x\in\R^d}\int_{\R^d}{g^2(x)\over\vert y-x\vert^p}\,dy.\nonumber
\end{eqnarray}
Clearly, $\rho^*_{Q_{r+1}}\le\rho^*(d,p)$. By
(\ref{poisson-19}), $\rho^*(d,p)=\rho(d,p)$.\vspace*{1pt}

By (\ref{theorem-7'}) in Lemma \ref{A-1}, therefore,
\begin{eqnarray}\label{u-12}
&&\lim_{t\to\infty}{1\over\log t}\log
\P\biggl\{\sup_{g\in\mathcal{G}_d(Q_{r+1})}\int_{Q_{r+1}}\xi_t(x)g^2(x)\,dx
>1 \biggr\}\nonumber\\[-8pt]\\[-8pt]
&& \qquad\le-p^{-p/2}(2-p)^{-{(4-p)/2}} \biggl({2u\over1+\delta}
\biggr)^{(2-p)/2}{2+d-p\over\theta\sigma(d,p)}.\nonumber
\end{eqnarray}

Again, (\ref{u-11}) [in the context of (\ref{u-2})]
follows forms (\ref{u-7}), (\ref{u-8}) and (\ref{u-12}).\vadjust{\goodbreak}

For any $\gamma>1$, (\ref{u-11}) implies that
\[
\sum_{k}\P\bigl\{\lambda_{\theta\xi}(Q_{\gamma^k})\ge\bigl(\Lambda(\theta
)+\delta\bigr)
h_{\gamma^k}^2 \bigr\}<\infty.
\]
By the Borel--Cantelli lemma,
\[
\limsup_{k\to\infty}h_{\gamma^k}^{-2}\lambda_{\theta\xi}(Q_{\gamma
^k})\le
(1+2\delta)\Lambda(\theta)\qquad  \mbox{a.s.}
\]
Since $\lambda_{\theta\xi}(Q_{t})$ is monotonic in $t$ and
$\delta>0$ can be arbitrarily small, we have proved (\ref{u-6}).

\section{Lower bounds}\label{l}

In this section we prove that
%
\begin{eqnarray}\label{l-1}
\liminf_{t\to\infty}t^{-1}(\log t)^{-{(d-p)/ d}}\log
\E_0\exp\biggl\{-\theta\int_0^t
\ol{V}(B_s)\,ds \biggr\}\ge\Lambda_0(\theta)  \nonumber\\[-9pt]\\[-9pt]
\eqntext{\mbox{a.s.-}\P}
\end{eqnarray}
when $d/2<p<d$ and
%
\begin{eqnarray}\label{l-2}
\liminf_{t\to\infty}{1\over t} \biggl({\log\log t\over\log t} \biggr)^{2/ (2-p)}
\log\E_0\exp\biggl\{\theta\int_0^t\ol{V}(B_s)\,ds \biggr\}
\ge\Lambda_1(\theta)\nonumber\\[-9pt]\\[-9pt]
\eqntext{\mbox{a.s.-}\P}
\end{eqnarray}
when $d/2<p<\min\{2, d\}$;
where $\Lambda_0(\theta)$ and $\Lambda_1(\theta)$
are given in~(\ref{u-1'}) and~(\ref{u-2'}), respectively.

Let $h_t$ be defined in (\ref{u-3}), and write $\xi(x)=-\ol{V}(x)$
in connection with the proof of
(\ref{l-1}) and $\xi(x)=\ol{V}(x)$ in connection with the proof of
(\ref{l-2}). Let $0<q<1$ be fixed but close to 1. Let $\alpha,\beta>1$
satisfy $\alpha^{-1}+\beta^{-1}=1$ with $\alpha$ being close to 1.
By (\ref{FK-5}) in Lemma \ref{FK-3},
%
\begin{eqnarray} \label{l-4}
&&\E_0\exp\biggl\{\theta\int_0^t\xi(B_s)\,ds \biggr\}\nonumber\\[-2pt]
&&\qquad\ge\biggl(\E_0\exp\biggl\{-{\theta\beta\over\alpha}\int_0^{t^q}\xi(B_s)\,
ds \biggr\}
\biggr)^{-\alpha/\beta}\nonumber\\[-2pt]
&&\qquad\quad{}\times\biggl\{\int_{Q_{t^q}}p_{t^q}(x)\E_x \biggl[\exp\biggl\{\alpha^{-1}
\int_0^{t-t^q}\xi(B_s)\biggr\}; \tau_{Q_{t^q}}\ge t-t^q \biggr]\,dx \biggr\}^\alpha
\nonumber\\[-2pt]
&&\qquad\ge{1\over(2\pi t^q)^{\alpha d/2}}e^{-c t^q}
\biggl(\E_0\exp\biggl\{-{\theta\beta\over\alpha}\int_0^{t^q}\xi(B_s)\,ds \biggr\}
\biggr)^{-\alpha/\beta}\\[-2pt]
&&\qquad\quad{}\times\biggl\{\int_{Q_{t^q}}\E_x \biggl[\exp\biggl\{\alpha^{-1}
\int_0^{t-t^q}\xi(B_s)\biggr\}; \tau_{Q_{t^q}}\ge t-t^q \biggr]\,dx \biggr\}^\alpha
\nonumber\\[-2pt]
&&\qquad\ge e^{-c_1 t^q}
\biggl(\E_0\exp\biggl\{-{\theta\beta\over\alpha}\int_0^{t^q}\xi(B_s)\,ds \biggr\}
\biggr)^{-\alpha/\beta}\nonumber\\[-2pt]
&&\qquad\quad{}\times\exp\bigl\{-(\alpha^2/\beta)t^q\lambda_{(\beta/\alpha
^2)\theta\xi}(Q_{t^q})
+\alpha^2 t\lambda_{\alpha^{-2}\theta\xi}(Q_{t^q}) \bigr\}\nonumber[-2pt]
\end{eqnarray}
for large $t$,
where the last step follows from Lemma \ref{FK-1'} (with $\delta=t^q$
and $t$ being replaced by $t-t^q$), and
the positive constant $c_1$ is made to be larger than $c$ for absorbing all
bounded-by-polynomial quantities including those
appearing on the right-hand side of
(\ref{FK-2'}).

By (\ref{u-1}), (\ref{u-2}) and (\ref{u-6})
\begin{eqnarray*}
\log\E_0\exp\biggl\{-{\theta\beta\over\alpha}\int_0^{t^q}\xi(B_s)\,ds \biggr\}
&=&o(t) \quad\mbox{and}\\
\lambda_{(\beta/\alpha^2)\theta\xi}(Q_{t^q})
&=&O(h_t^2) \qquad\mbox{a.s.}
\end{eqnarray*}
as $t\to\infty$. Therefore, all we need is to show that
%
\begin{equation}\label{l-5}
\liminf_{t\to\infty}h_t^{-2}\lambda_{\theta\xi}(Q_t)\ge\Lambda(\theta)
\qquad\mbox{a.s.}
\end{equation}
for every $\theta>0$, where $\Lambda(\theta)$ is given in (\ref{u-6}).
Indeed, applying (\ref{l-5}) to (\ref{l-4}) with $\theta$ being
replaced by
$\alpha^{-2}\theta$ leads to
\[
\liminf_{t\to\infty}t^{-1}h_{t^q}^{-2}
\log\E_0\exp\biggl\{\theta\int_0^t\xi(B_s)\,ds \biggr\}
\ge\alpha^2\Lambda(\alpha^{-2}\theta) \qquad\mbox{a.s.}
\]
Letting $\alpha\to1^+$, the right-hand side tends to $\Lambda(\theta)$.
In addition, $h_{t^q}=q^{(d-p)/(2d)}h_t$ when applied to (\ref{l-1})
and $h_{t^q}\sim q^{1/(2-p)}h_t$ when applied to (\ref{l-2}). Therefore,
with probability 1,
\begin{eqnarray*}
&&\liminf_{t\to\infty}t^{-1}h_{t}^{-2}
\log\E_0\exp\biggl\{\theta\int_0^t\xi(B_s)\,ds \biggr\}\\
&&\qquad \ge\cases{
q^{(d-p)/d}\Lambda_0(\theta),&\quad
when applied to (\ref{l-1}),\cr
q^{2/(2-p)}\Lambda_1(\theta),&\quad
when applied to (\ref{l-2}).
}
\end{eqnarray*}
Letting $q\to1^-$ on the right-hand side leads to (\ref{l-1}) and (\ref{l-2}).

We now prove (\ref{l-5}). Let $u>0$ be fixed but arbitrary. Write
$\hat{h}_t=\sqrt{u}h_t$ and
\[
\eta_t(x)=\pm\theta\hat{h}_t^{p-2}\int_{\R^d}{1\over\vert y-x\vert^p} [
\omega(\hat{h}^{-d}_t\,dy)-\hat{h}^{-d}_t\,dx ],
\]
where ``$-$'' is for the proof of (\ref{l-1}) and ``$+$'' is for the
proof of
(\ref{l-2}).
Under the substitution $g(x)\mapsto\hat{h}_t^{d/2}g(h_tx)$,
\begin{eqnarray*}
\lambda_{\theta\xi}(Q_t)=\hat{h}_t^2\lambda_{\eta_t}(Q_{t\hat{h}_t}).
\end{eqnarray*}
Consequently,
%
\begin{eqnarray}\label{l-6}
\P\{\lambda_{\theta\xi}(Q_t)\le uh_t^2 \}
&=&\P\{\lambda_{\eta_t}(Q_{t\hat{h}_t})\le1 \}\nonumber\\[-8pt]\\[-8pt]
&=&\P\biggl\{\sup_{g\in\mathcal{G}_d(Q_{t\hat{h}_t})}\int_{\R^d}
\eta_t(x)g^2(x)\,dx\le1 \biggr\},\nonumber
\end{eqnarray}
where the last step follows from Lemma \ref{A-1}.\vadjust{\goodbreak}

Let $s>{2+d-p\over d-p}$ and $r>0$\vspace*{-2pt}
be fixed. When $t$ is large, $z+Q_r\subset Q_{t\hat{h}_t}$
for each $z\in h_t^s\Z^d\cap Q_{t\hat{h}_t-r}$. Hence,
\begin{eqnarray*}
&&\sup_{g\in\mathcal{G}_d(Q_{t\hat{h}_t})}\int_{\R^d}
\eta_t(x)g^2(x)\,dx\\
&&\qquad \ge\sup_{g\in\mathcal{G}_d(z+Q_r)}\int_{\R^d}
\eta_t(x)g^2(x)\,dx,\qquad  z\in h_t^s\Z^d\cap Q_{t\hat{h}_t-r}.
\end{eqnarray*}
Thus
%
\begin{equation}\label{l-7}
\sup_{g\in\mathcal{G}_d(Q_{t\hat{h}_t})}\int_{\R^d}
\eta_t(x)g^2(x)\,dx
\ge\max_{z\in h_t^s\Z^d\cap Q_{t\hat{h}_t-r}}
\sup_{g\in\mathcal{G}_d(z+Q_r)}\int_{\R^d}
\eta_t(x)g^2(x)\,dx.\hspace*{-31pt}
\end{equation}

Let the smooth function $\alpha(\cdot)$: $\R^+\longrightarrow[0,1]$
be given as in Section \ref{poisson}. Given $a>0$, write
\[
K_{t,a}(x)= \cases{
\displaystyle {\alpha(a^{-1}(\hat{h}_t)^{(2+d-p)/ (d-p)}\vert x\vert)
\over\vert x\vert^p},&\quad
when applied to (\ref{l-1}),\vspace*{2pt}\cr
\displaystyle {\alpha (a^{-1}(d\log\hat{h}_t)^{1/p}
\vert x\vert)\over\vert x\vert^p},&\quad
when applied to (\ref{l-2}),
}
\]
and
\[
L_{t,a}(x)= \cases{
\displaystyle {1-\alpha(a^{-1}(\hat{h}_t)^{(2+d-p)/(d-p)}\vert x\vert)
\over\vert x\vert^p},&\quad
when applied to (\ref{l-1}),\cr
\displaystyle {1-\alpha (a^{-1}(d\log\hat{h}_t)^{1/p}
\vert x\vert)\over\vert x\vert^p},&\quad
when applied to (\ref{l-2}).
}
\]
By the equality
\begin{eqnarray*}
&&\int_{\R^d}
\eta_t(x)g^2(x)\,dx\\
&&\qquad =(\pm\theta)\hat{h}_t^{p-2}\int_{\R^d} \biggl[\int
_{z+Q_r}K_{a,t}(y-x)g^2(y)\,dy \biggr]
[\omega(\hat{h}_t^{-d}\,dx)-\hat{h}_t^{-d}\,dx ]\\
&&\qquad \quad {}+(\pm\theta)\hat{h}_t^{p-2}\int_{\R^d} \biggl[\int
_{z+Q_r}L_{a,t}(y-x)g^2(y)\,dy \biggr]
[\omega(\hat{h}_t^{-d}dx)-\hat{h}_t^{-d}\,dx ]\\
&&\qquad =\theta\hat{h}_t^{p-2} \bigl(A_z(g)+B_z(g) \bigr)\qquad  (\mbox{say})
\end{eqnarray*}
and by triangular inequality, the right-hand side of (\ref{l-7})
is no less than
\[
\hat{h}_t^{p-2} \Bigl\{
\max_{z\in h_t^s\Z^d\cap Q_{\hat{h}_t-r}}\sup_{g\in\mathcal{G}_d(z+Q_r)}A_z(g)
-\max_{z\in h_t^s\Z^d\cap Q_{\hat{h}_t-r}}\sup_{g\in\mathcal
{G}_d(z+Q_r)}\vert B_z(g)\vert\Bigr\}.
\]

In addition, the random variables
\[
\sup_{g\in\mathcal{G}_d(z+Q_r)}\vert B_z(g)\vert;\qquad
z\in h_t^s\Z^d\cap Q_{t\hat{h}_t-r},
\]
are identically distributed. Therefore, for any $\delta>0$,
\begin{eqnarray*}
&&\P\Bigl\{\max_{z\in h_t^s\Z^d\cap Q_{t\hat{h}_t-r}}
\sup_{g\in\mathcal{G}_d(z+Q_r)}\vert B_z(g)\vert\ge\delta\theta^{-1}
\hat{h}_t^{2-p} \Bigr\}\\
&&\qquad \le\# \{h_t^s\Z^d\cap Q_{t\hat{h}_t-r} \}
\P\Bigl\{\sup_{g\in\mathcal{G}_d(Q_r)}\vert B_0(g)\vert\ge
\delta\theta^{-1}\hat{h}_t^{2-p} \Bigr\}.
\end{eqnarray*}

Further, since $\alpha(\cdot)$ is supported on $[0,3]$ and $s>{2\over d-p}$,
\[
A_z(g)=\pm\theta\int_{z+Q_{2^{-1}h_t^s}} \biggl[\int
_{z+Q_r}K_{a,t}(y-x)g^2(y)\,dy \biggr]
[\omega(\hat{h}_t^{-d}\,dx)-\hat{h}_t^{-d}\,dx ]
\]
for all $z\in h_t^s\Z^d\cap Q_{t\hat{h}_t-r}$
as $t$ is sufficiently large. Consequently,
the random variables
\[
\sup_{g\in\mathcal{G}_d(z+Q_r)}A_z(g);\qquad
z\in h_t^s\Z^d\cap Q_{t\hat{h}_t-r},
\]
form an i.i.d. sequence. Therefore,
\begin{eqnarray*}
&&\P\biggl\{\max_{z\in h_t^s\Z^d\cap Q_{t\hat{h}_t-r}}\sup_{g\in\mathcal{G}_d(z+Q_r)}
A_z(g)\le{1+\delta\over\theta}\hat{h}_t^{2-p} \biggr\}\\
&&\qquad = \biggl(\P\biggl\{\sup_{g\in\mathcal{G}_d(Q_r)}A_0(g)\le
{1+\delta\over\theta}
\hat{h}_t^{2-p} \biggr\}
\biggr)^{\# \{h_t^s\Z^d\cap Q_{t\hat{h}_t-r} \}}\\
&&\qquad = \biggl(1-\P\biggl\{\sup_{g\in\mathcal{G}_d(Q_r)}A_0(g)\ge
{1+\delta\over\theta}\hat{h}_t^{2-p} \biggr\}
\biggr)^{\# \{h_t^s\Z^d\cap Q_{t\hat{h}_t-r} \}}.
\end{eqnarray*}

Summarizing our argument since (\ref{l-6}) and (\ref{l-7}),
%
\begin{eqnarray}\label{l-8}
&&\P\{\lambda_{\theta\xi}(Q_t)\le uh_t^2 \}\nonumber\\
&&\qquad\le \biggl(1-\P\biggl\{\sup_{g\in\mathcal{G}_d(Q_r)}A_0(g)\ge
{1+\delta\over\theta}\hat{h}_t^{2-p} \biggr\}
\biggr)^{\# \{h_t^s\Z^d\cap Q_{t\hat{h}_t-r} \}}\\
&&\quad\qquad{}+\# \{h_t^s\Z^d\cap Q_{t\hat{h}_t-r} \}
\P\Bigl\{\sup_{g\in\mathcal{G}_d(Q_r)}\vert B_0(g)\vert\ge
\delta\theta^{-1}\hat{h}_t^{2-p} \Bigr\}.\nonumber
\end{eqnarray}

Once again, we reach the point of using Theorem \ref{poisson-6}
and Theorem \ref{poisson-14}. In connection with (\ref{l-1}), by definition
\begin{eqnarray*}
\sup_{g\in\mathcal{G}_d(Q_r)}A_0(g)&=&-\inf_{g\in\mathcal{G}_d(Q_r)}
\int_{\R^d} \biggl[\int_{Q_r}K_{a,t}(y-x)g^2(y)\,dy \biggr] \\
&&{}\times[\omega(\hat{h}_t\,dx)-\hat{h}_t^{-d}\,dx ],
\\
\sup_{g\in\mathcal{G}_d(Q_r)}\vert B_0(g)\vert&=&\sup_{g\in\mathcal{G}_d(Q_r)}
\biggl\vert\int_{\R^d} \biggl[\int_{Q_r}L_{a,t}(y-x)g^2(y)\,dy \biggr]\\
&&\hspace*{61pt}{}\times[\omega(\hat{h}_t^{-d}\,dx)-\hat{h}_t^{-d}\,dx ] \biggr\vert.
\end{eqnarray*}
Taking $\epsilon=\hat{h}_t^{-d}$ in (\ref{poisson-8}) and (\ref{poisson-7}),
\begin{eqnarray*}
&&\liminf_{a\to\infty}\liminf_{t\to\infty}{1\over\log t}\log
\P\biggl\{\sup_{g\in\mathcal{G}_d(Q_r)}A_0(g)\ge
{1+\delta\over\theta}\hat{h}_t^{2-p} \biggr\}\\
&&\qquad \ge-u^{d/(d-p)}I_{Q_r} \biggl({1+\delta\over\theta} \biggr),
\\
&&\lim_{a\to\infty}\limsup_{t\to\infty}{1\over\log t}
\log\P\Bigl\{\sup_{g\in\mathcal{G}_d(Q_r)}\vert B_0(g)
\vert\ge\delta\theta^{-1}\hat{h}_t^{2-p} \Bigr\}=-\infty,
\end{eqnarray*}
where the rate functions $I_{Q_r}(\cdot)$ are defined in (\ref{poisson-10}).

By definition,
\[
\sup_{g\in\mathcal{G}_d(Q_r)}
\|g\|_{\mathcal{L}^2(Q_r)}^2\le1.
\]
We claim that
%
\begin{equation}\label{l-9}
\lim_{r\to\infty}\sup_{g\in\mathcal{G}_d(Q_r)}
\|g\|_{\mathcal{L}^2(Q_r)}^2=1.
\end{equation}
Indeed, for a fixed $g\in\mathcal{F}_d(Q_{1})$ the function
\[
f_r(x)= \biggl(r^d+{1\over2}r^{d-2}\|\nabla g\|_{\mathcal{L}^2(Q_1)}^2 \biggr)^{-1/2}
g \biggl({x\over r} \biggr),\qquad  x\in Q_r,
\]
is in $\mathcal{G}_d(Q_r)$ and
\begin{eqnarray}
\sup_{g\in\mathcal{G}_d(Q_r)}
\|g\|_{\mathcal{L}^2(Q_r)}^2\ge\|f_r\|_{\mathcal{L}^2(Q_{r})}^2
={r^d\over r^d+{(1/2)}r^{d-2}\|\nabla g\|_{\mathcal{L}^2(Q_1)}^2}
\longrightarrow1\nonumber \\
\eqntext{(r\to\infty).}
\end{eqnarray}

By (\ref{l-9}) and by the definition of $I_{Q_r}(\cdot)$ given in
(\ref{poisson-10}),
\[
\lim_{r\to\infty} I_{Q_r}
\biggl({1+\delta\over\theta} \biggr)
= \biggl({(d-p)(1+\delta)\over d\theta} \biggr)^{d/(d-p)}
\biggl({\omega_d\over d}\Gamma\biggl({2p-d\over p} \biggr)
\biggr)^{-{p/( d-p)}}.
\]

Take $u<(1+2\delta)^{-1}\Lambda_0(\theta)$.
There is a $\nu(\delta)>0$ such that when $a$ and $r$ are sufficiently large,
\[
\P\biggl\{\sup_{g\in\mathcal{G}_d(Q_r)}A_0(g)\ge
{1+\delta\over\theta}\hat{h}_t^{2-p} \biggr\}
\ge\exp\bigl\{- \bigl(d-\nu(\delta) \bigr)\log t \bigr\}
=t^{- (d-\nu(\delta) )}
\]
and
\[
\P\Bigl\{\sup_{g\in\mathcal{G}_d(Q_r)}\vert B_0(g)
\vert\ge\delta\theta^{-1}\hat{h}_t^{2-p} \Bigr\}\le\exp\{-2d\log t\}
=t^{-2d}
\]
for sufficiently large $t$.

Being brought to (\ref{l-8}), our estimates give
%
\begin{eqnarray}\label{l-9'}
&&\P\{\lambda_{\theta\xi}(Q_t)\le uh_t^2 \}\nonumber\hspace*{-35pt}\\
&&\qquad\le \bigl(1-t^{-(d-\nu(\delta))}
\bigr)^{\# \{h_t^s\Z^d\cap Q_{t\hat{h}_t-r} \}}
+\# \{h_t^s\Z^d\cap Q_{\sqrt{u}th_t-r} \}
t^{-2d}\hspace*{-35pt}\\
&&\qquad\le\exp\bigl\{-c_1t^{\nu(\delta)}h_t^{-d(s-1)} \bigr\}
+c_2 t^{-d}.\nonumber\hspace*{-35pt}
\end{eqnarray}

For any $\gamma>1$ and $u<(1+2\delta)^{-1}\Lambda_0(\theta)$, therefore,
\[
\sum_k\P\{\lambda_{\theta\xi}(Q_{\gamma^k})\le uh_{\gamma^k}^2 \}
<\infty.
\]
By the Borel--Cantelli lemma,
\[
\liminf_{k\to\infty}h_{\gamma^k}^{-2}\lambda_{\theta\xi}(Q_{\gamma
^k})\ge
(1+2\delta)^{-1}\Lambda_0(\theta) \qquad\mbox{a.s.}
\]
Since $\lambda_{\theta\xi}(Q_t)$ is monotonic in $t$ and $\delta>0$ can be
arbitrarily small, we have proved (\ref{l-5}) associated with (\ref{l-1}).

As for (\ref{l-2}), by (\ref{poisson-15}) and
(\ref{poisson-16}) (with $\epsilon=\hat{h}_t^{-d}$),
\begin{eqnarray*}
&&\lim_{t\to\infty}{1\over\log t}\log\P\Bigl\{\sup_{g\in\mathcal
{G}_d(Q_r)}\vert B_0(g)
\vert\ge\delta\theta^{-1}\hat{h}_t^{2-p} \Bigr\}
=-\infty,
\\
&&\lim_{t\to\infty}{1\over\log t}\log\P\Bigl\{\sup_{g\in\mathcal{G}_d(Q_r)} A_0(g)
\ge{1+\delta\over\theta}\hat{h}_t^{2-p} \Bigr\}\\
&&\qquad =-u^{d/(d-p)}{2+d-p\over(2-p)\rho^*_{Q_r}}{1+\delta\over\theta}
\ge-u^{d/(d-p)}{2+d-p\over(2-p)\rho_{Q_r}}{1+\delta\over\theta},
\end{eqnarray*}
where $\rho_D^*$ is defined in (\ref{poisson-13}) and
\[
\rho_{Q_r}=\sup_{g\in\mathcal{G}_d(Q_r)}\int_{Q_r}{g^2(x)\over\vert
x\vert^p}\,dx.
\]
Clearly, $\rho_{Q_r}$ is nondecreasing in $r$ and $\rho_{Q_r}\le\rho(d,p)$,
where $\rho(d,p)$ is defined in (\ref{u-11''}).
We claim that
%
\begin{equation}\label{l-12}
\lim_{r\to\infty}\rho_{Q_r}=\rho(d,p).
\end{equation}
Indeed, let $\alpha(\cdot)$ be the smooth truncation function
introduced in Section \ref{poisson}. For any $f\in\mathcal{G}_d(\R^d)$,
write
\[
f_r(x)=f(x)\alpha(3r^{-1}\vert x\vert).
\]
The function
%
\begin{equation}\label{l-13}
g_r(x)= \bigl(\|f_r\|_{\mathcal{L}^2(Q_r)}+2^{-1}
\|\nabla f_r\|_{\mathcal{L}^2(Q_r)} \bigr)^{-1/2}f_r(x)
\end{equation}
is in $\mathcal{G}_d(Q_r)$. Thus, by the fact that
$\alpha(\cdot)\ge1_{[0,1]}(\cdot)$
\[
\rho_{Q_r}\ge\int_{Q_r}{g_r^2(x)\over\vert x\vert^p}\,dx
\ge\biggl(\|f_r\|_{\mathcal{L}^2(Q_r)}^2+{1\over2}
\|\nabla f_r\|_{\mathcal{L}^2(Q_r)}^2 \biggr)^{-1}\int_{\{\vert x\vert\le
r/3\}}
{f^2(x)\over\vert x\vert^p}\,dx.
\]
Notice that
$\|f_r\|^2_{\mathcal{L}^2(Q_r)}\le\|f\|_2^2$ and
\begin{eqnarray*}
\vert\nabla f_r(x)\vert
&\le&3r^{-1}\vert\alpha'(3r^{-1}\vert x\vert)\vert\cdot\vert f(x)\vert
+\alpha(3r^{-1}\vert x\vert)\vert\nabla f(x)\vert\\
&\le& 3r^{-1}\vert f(x)\vert+\vert\nabla f(x)\vert,
\end{eqnarray*}
where the last step follows from the fact that
$\vert\alpha(\cdot)\vert\le1$ and $\vert\alpha'(\cdot)\vert\le1$.

Thus,
%
\begin{equation}\label{l-14}
\qquad \liminf_{r\to\infty} \bigl(\|f_r\|^2_{\mathcal{L}^2(Q_r)}+\tfrac{1}{2}
\|\nabla f_r\|^2_{\mathcal{L}^2(Q_r)} \bigr)^{-1}\ge\bigl(\|f\|_2^2+\tfrac{1}{2}
\|\nabla f\|_2^2 \bigr)^{-1}=1.
\end{equation}

Summarizing our argument,
\[
\liminf_{r\to\infty}\rho_{Q_r}\ge\int_{\R^d}
{f^2(x)\over\vert x\vert^p}\,dx.
\]
Taking supremum over $f\in\mathcal{G}_d$ on the right-hand side
leads to (\ref{l-12}).

By (\ref{l-12}) and (\ref{theorem-7'}), therefore,
\[
\lim_{r\to\infty}\rho_{Q_r}= \biggl({2-p\over2} \biggr)^{(2-p)/2}p^{p/2}\sigma(d,p).
\]

Similarly, the above discussion leads
to (\ref{l-5}) [corresponding to (\ref{l-2})], again by the
Borel--Cantelli lemma.

\begin{appendix}
\section*{Appendix}\label{A}
\setcounter{equation}{0}

\begin{lemma}\label{psi-1}
Under $d/2<p<d$,
%
\begin{equation}\label{psi-2}
\int_{\R^d} \biggl[\exp\biggl\{-{1\over\vert x\vert^p} \biggr\}
-1+{1\over\vert x\vert^p} \biggr]\,dx=\omega_d{p\over d-p}
\Gamma\biggl({2p-d\over p} \biggr),
\end{equation}
where $\omega_d$ is the volume of the $d$-dimensional unit ball.
\end{lemma}

\begin{pf}
By the sphere substitution,
\begin{eqnarray*}
\int_{\R^d} \biggl[\exp\biggl\{-{1\over\vert x\vert^p} \biggr\}
-1+{1\over\vert x\vert^p} \biggr]\,dx&=&d\omega_d\int_0^\infty
\biggl[\exp\biggl\{-{1\over\rho^p} \biggr\}
-1+{1\over\rho^p} \biggr]\rho^{d-1}\,d\rho\\
&=&{d\omega_d\over p}\int_0^\infty [e^{-\gamma}-1+\gamma]
\gamma^{-{(d+p)/p}}\,d\gamma,
\end{eqnarray*}
where the second step follows from the substitution
$\rho=\gamma^{-1/p}$.\vadjust{\goodbreak}

Applying the integration by parts twice (under the assumption \mbox{$d/2\,{<}\,p\,{<}\,d$}),
\begin{eqnarray*}
\int_0^\infty [e^{-\gamma}-1+\gamma]
\gamma^{-{(d+p)/p}}\,d\gamma&=&{p\over d}\int_0^\infty[1-e^{-\gamma
} ]
\gamma^{-d/p}\,d\gamma\\
&=&{p^2\over d(d-p)}\int_0^\infty\gamma^{-{(d-p)/d}}e^{-\gamma
}\,d\gamma\\
&=&{p^2\over d(d-p)}\Gamma\biggl({2p-d\over p} \biggr).
\end{eqnarray*}
We have proved identity (\ref{psi-2}).
\end{pf}

Recall that for any domain $D\subset\R^d$,
\begin{eqnarray*}
\mathcal{G}_d(D)&=& \bigl\{g\in W^{1,2}(D); \|g\|_{\mathcal{L}^2(D)}^{2}+
\tfrac{1}{2}\|\nabla g\|_{\mathcal{L}^2(D)}^{-2}=1 \bigr\},
\\
\mathcal{F}_d(D)&=& \bigl\{g\in W^{1,2}(D); \|g\|_{\mathcal{L}^2(D)}=1 \bigr\}.
\end{eqnarray*}
In particular, $\mathcal{G}_d=\mathcal{G}_d(\R^d)$ and $\mathcal
{F}_d=\mathcal{F}_d(\R^d)$.

\begin{lemma}\label{A-1} Let the functional $Z(g^2)$ [$g\in
W^{1,2}(D)$]
satisfy
$Z(cg^2)=cZ(g^2)$ for every $g\in W^{1,2}(D)$ and $c>0$. Then
\[
\sup_{g\in\mathcal{F}_d(D)} \biggl\{Z(g^2)-{1\over2}\int_D\vert\nabla
g(x)\vert^2\,dx\biggr\}>1
\]
if any only if $ \sup_{g\in\mathcal{G}_d(D)}Z(g^2)>1$.
\end{lemma}

\begin{pf}
For any $g\in\mathcal{F}_d(D)$,
\[
Z(g^2)
\le \Bigl(\sup_{f\in\mathcal{G}_d(D)}Z(f^2) \Bigr)
\biggl(1+{1\over2}\int_D\vert\nabla g(x)\vert^2\,dx \biggr).
\]
Hence,
\begin{eqnarray*}
\hspace*{-4pt}&&\sup_{g\in\mathcal{F}_d(D)} \biggl\{Z(g^2)
-{1\over2}\int_D\vert\nabla g(x)\vert^2\,dx
\biggr\}\\
\hspace*{-4pt}&&\qquad \le\sup_{g\in\mathcal{F}_d(D)} \biggl\{ \biggl(\sup_{f\in\mathcal{G}_d(D)}Z(f^2) \biggr)
\biggl(1+{1\over2}\int_D\vert\nabla g(x)\vert^2\,dx \biggr)-{1\over2}
\int_D\vert\nabla g(x)\vert^2\,dx \biggr\}.
\end{eqnarray*}
Therefore, $ \sup_{g\in\mathcal{G}_d(D)}Z(g^2)>1$, if
\[
\sup_{g\in\mathcal{F}_d(D)} \biggl\{Z(g^2)-{1\over2}\int_D\vert\nabla
g(x)\vert^2\,dx\biggr\}>1.
\]

On the other hand, assume $ \sup_{g\in\mathcal{G}_d(D)}Z(g^2)>1$.
Then there is $g_0\in\mathcal{G}_d(D)$ such that
$Z(g_0^2)>1$. Write
$f_0(x)=g_0(x)/\|g_0\|_{\mathcal{L}^2(D)}$. We
have $f_0\in\mathcal{F}_d(D)$\vadjust{\goodbreak} and
\[
Z(f_0^2)-{1\over2}\int_D\vert\nabla f_0(x)\vert^2\,dx
>\|g_0\|_{\mathcal{L}^2(D)}^{-2}-\|g_0\|_{\mathcal{L}^2(D)}^{-2}
\bigl(1-\|g_0\|_{\mathcal{L}^2(D)}^2 \bigr)=1.\qquad
\]
\upqed
\end{pf}

It was shown (see \cite{BCR}, (1.19)) that for every $\lambda>0$,
%
\begin{equation}\label{theorem-4}
M(\lambda)\equiv\sup_{g\in\mathcal{F}_d} \biggl\{\lambda\int_{\R^d}
{g^2(x)\over\vert x\vert^p}\,dx-{1\over2}\int_{\R^d}\vert\nabla g(x)\vert
^2\,dx\biggr\}<\infty
\end{equation}
under $d/2<p<\min\{2,d\}$.

Further, by rescaling $g(x)\mapsto a^{d/2}g(ax)$ for suitable $a>0$,
one can show that
%
\begin{equation}\label{theorem-5}
M(\lambda)=\lambda^{2/(2-p)}M(1).
\end{equation}

\begin{lemma}\label{th-1} Under $d/2<p<\min\{2,d\}$, there is a constant
$C>0$ such that
%
\begin{equation}\label{theorem-6}
\int_{\R^d}{f^2(x)\over\vert x\vert^p}\,dx
\le C\|f\|_2^{2-p}\|\nabla f\|_2^{p} \qquad\forall f\in W^{1,2}(\R^d).
\end{equation}
Further, let $\sigma(d,p)$ be the smallest (infimum) among above $C$. Then
%
\begin{equation}\label{theorem-7}
M(\lambda)={2-p\over2}p^{p/(2-p)}
(\lambda\sigma(d,p) )^{2/(2-p)},\qquad  \lambda>0.
\end{equation}

In addition,
%
\begin{equation}\label{theorem-6'}
\rho(d,p)\equiv\sup\biggl\{\int_{\R^d}{g^2(x)\over\vert x\vert^p}\,dx;
g\in\mathcal{G}_d \biggr\}<\infty
\end{equation}
and
%
\begin{equation}\label{theorem-7'}
\rho(d,p)= \biggl({2-p\over2} \biggr)^{(2-p)/2}p^{p/2}\sigma(d,p).
\end{equation}
\end{lemma}

\begin{pf}
In view of (\ref{theorem-5}) we may take $\lambda=1$
in (\ref{theorem-7}). For any $f\in W^{1,2}$ with $\|f\|_2=1$,
let
\[
\int_{\R^d}{f^2(x)\over\vert x\vert^p}\,dx
=C_f\|\nabla f\|_2^{p}.
\]
Given $\gamma>0$, let $g(x)=\gamma^{d/2}f(\gamma x)$. Then $\|g\|_2=1$,
$\|\nabla g\|_2=\gamma\|\nabla f\|_2$, and therefore
\[
\int_{\R^d}{g^2(x)\over\vert x\vert^p}\,dx=\gamma^p
\int_{\R^d}{f^2(x)\over\vert x\vert^p}\,dx=\gamma^pC_f\|\nabla f\|_2^{p}
=C_f\|\nabla g\|_2^{p}.
\]

Thus
\[
M(1)\ge C_f\|\nabla g\|_2^{p}
-\tfrac{1}{2}\|\nabla g\|_2^{2}=C_f\gamma^p\|\nabla f\|_2^{p}
-\tfrac{1}{2}\gamma^2\|\nabla f\|_2^{2}.
\]
Since $\gamma>0$ is arbitrary, the variable $\gamma\|\nabla f\|_2$
runs over all positive numbers. Consequently,
\[
M(1)\ge\sup_{x>0} \biggl\{ C_fx^{p}-{1\over2}
x^2 \biggr\}={2-p\over2}C_f^{2/(2-p)}p^{p/(2-p)}.
\]
By homogeneity, we have proved (\ref{theorem-6})
with
\[
M(1)\ge{2-p\over2}p^{p/(2-p)}\sigma(d,p)^{2/(2-p)}.
\]

On the other hand, for any $g\in\mathcal{F}_d$
\begin{eqnarray*}
\int_{\R^d}
{g^2(x)\over\vert x\vert^p}\,dx-{1\over2}\int_{\R^d}\vert\nabla g(x)\vert^2\,dx
&\le&\sigma(d,p)\|\nabla g\|_2^p-{1\over2}\|\nabla g\|_2^2\\
&\le&\sup_{x>0}
\biggl\{\sigma_1(d,p)x^p-{1\over2}x^2 \biggr\}\\
&=&{2-p\over2}p^{p/(2-p)}\sigma(d,p)^{2/(2-p)}.
\end{eqnarray*}
We have proved (\ref{theorem-7}).

Obviously, (\ref{theorem-6'}) follows from (\ref{theorem-6}). Take
\[
Z(g^2)={1\over\rho(d,p)}\int_{\R^d}{g^2(x)\over\vert x\vert^p}\,dx,\qquad
g\in W^{1,2}(\R^d).
\]
We have that $ \sup_{g\in\mathcal{G}_d}Z(g^2)=1$.
By (\ref{theorem-5}), the function $M(\lambda)$
is continuous and increasing. By Lemma \ref{A-1}, we must have
%
\begin{equation}\label{theorem-10'}
M \biggl({1\over\rho(d,p)} \biggr)=1.
\end{equation}

Finally, (\ref{theorem-7'}) follows from (\ref{theorem-7}) and
(\ref{theorem-10'}).
\end{pf}

Another variation appearing in this paper is
\[
\rho^*(d,p)=\sup_{g\in\mathcal{G}_d}\sup_{x\in\R^d}
\int_{\R^d}{g^2(y)\over\vert y-x\vert^p}\,dy.
\]
We now claim that
%
\begin{equation}\label{poisson-19}
\rho^*(d,p)=\rho(d,p).
\end{equation}
Indeed,
\begin{eqnarray*}
\rho^*(d,p)&=&\sup_{x\in\R^d}\sup_{g\in\mathcal{G}_d}
\int_{\R^d}{g^2(y)\over\vert y-x\vert^p}\,dy=\sup_{x\in\R^d}\sup_{g\in
\mathcal{G}_d}
\int_{\R^d}{g_x^2(y)\over\vert y\vert^p}\,dy\\
&\le&\sup_{g\in\mathcal{G}_d}
\int_{\R^d}{g^2(y)\over\vert y\vert^p}\,dy=\rho(d,p),
\end{eqnarray*}
where $g_x(y)=g(x+y)$, and the inequality follows from the fact that
$g_x\in\mathcal{G}_d$ as soon as $g\in\mathcal{G}_d$.
\end{appendix}

\section*{Acknowledgment}
The author would like to thank
the anonymous referee who read the first version
of this paper for his/her interesting remarks and
excellent suggestions.


%

\printaddresses

\end{document}